\documentstyle[amsmath,amssymb]{article}

\makeatletter
\@addtoreset{equation}{subsection}
\makeatother

\begin{document}

\title{Conditional expectations on compact quantum groups and new examples
of quantum hypergroups.}

\author{A.A.Kalyuzhnyi$^*$}
\date{}
\maketitle
\begin{abstract}
\par We propose a construction of quantum hypergroups using
conditional expectations on compact quantum groups.
Using this construction, we describe several series of non-trivial
finite-dimensional quantum hypergroups via conditional expectations
of Delsart type on  non-trivial Kac algebras obtained by twisting of
the finite groups.
\end{abstract}
\footnotetext[1]
         {$^*$ This work has been partially supported by the
                Ukrainian Committee for Fundamental studies
                }

\subsection*{0. Introduction}
\setcounter{subsection}{0}
\setcounter{equation}{0}
The compact quantum hypergroup was introduced in \cite{ChaVai} as a structure
that simultaniosly generalies usual hypergroups (\cite{BK}, \cite{BH}),
compact quantum groups (\cite{Wor})
and bialgebras of biinvariant functions associated with
quantum Gelfand pairs (\cite{Koo}, \cite{ChaVai:double}, \cite{Vai:double},
\cite{Vai:double1}).
A compact quantum hypergroup is a unital $C^*$-algebra equipped with
a coassociative completely positive coproduct that preserves the unit element
and satisfies some additional axioms.
In \cite{ChaVai} a theory of corepresentations was studied and
there was established an analog of the Peter-Weyl theory
for compact quantum hypergroups.

Howewer, as we know, the only {\it concrete} nontrivial example of the
compact quantum hypergroup was constucted in \cite{PVai}.
The main purpose of this paper is to build a sufficient number
of such examples. We propose a general construction of
quantum hypergroups using conditional expectations on compact quantum groups
and, using this construction, we carry, to the quantum case, the Delsart construction
of hypergroups (recall that if $G$ is a locally compact group
and $\Gamma$ is a compact group of automorphisms of $G$ with a Haar measure
$\mu_\gamma$, then the algebra of continuous
$\Gamma$-invariant functions on $G$ equipped with the coproduct
$$
  (\Delta f)(x, y) = \int_{\Gamma} f(x\gamma (y)) d\mu_\gamma,
$$
is a hypergroup).

Using the Delsart construction we build several series of non-trivial
finite-di\-men\-sional quantum hypergroups from  non-trivial Kac algebras
obtained in \cite{Vai:twist}, \cite{Nik:twist}
by twisting the classical series of finite groups.

The paper is organized as follows. Section 1 contains neáessary
preliminaries on compact quantum hypergroups and on twisting
the finite groups.
In Section 2 we describe
the construction of compact quantum hypergroups using conditional expectations
on quantum hypergroups (and, in particular, on  compact quantum groups).
In the commutative case, we compare this construction
with orbital morphisms \cite{Jew}. In Section 3
we show that this construction includes the double coset construction in
quantum groups \cite{ChaVai:double}, \cite{Vai:double},
\cite{Vai:double1} and an analogue of the Delsart
construction for quantum hypergroup.
We also give here sufficient conditions
for an authomorphism of the initial group (or a certain automorphism
constructed from it) to be an automorphism of the corresponding
twisted Kac algebra. Section 4 contains five series of nontrivial
examples of finite-dimensional Delsart quantum hypergroups associated
with twistings of the classical series of finite groups.

I am deeply grateful to L. I. Vainerman and Yu. A. Chapovsky
for many stimulating discussions.

\subsection*{1. Preliminaries}
\setcounter{subsection}{1}
\setcounter{equation}{0}
{\bf 1.1.} Let $(A,\cdot ,1,*)$ be a separable unital $C^*$-algebra. We
denote by $A\otimes A$ the injective or projective $C^*$-tensor
square of $A$. We will call $(A,\Delta ,\epsilon,\star)$ a {\em
hypergroup structure} on the $C^*$-algebra $(A,\cdot ,1,*)$
if:
\begin{itemize}
  \item [(HS$_1$)]
               $(A,\Delta,\epsilon,\star)$ is a $\star$-coalgebra
               with a counit $\epsilon$, i.e. $\Delta:A\to A\otimes
		A$ and $\epsilon:A\to {\bf C}$ are linear mappings,
		$\star:A\to A$ is an antilinear mapping such that
		\begin{eqnarray}
			&
                        (\Delta\otimes id)\circ \Delta=
                        (id\otimes \Delta)\circ \Delta,&
			\label{co1}\\
                        &(\epsilon\otimes id)\circ \Delta=
                        (id\otimes \epsilon)\circ \Delta=id,&
			\label{co2}\\
                        &\Delta\circ \star =\Sigma\circ ( \star \otimes
                        \star )\circ \Delta,&\label{co3}\\
			&\star \circ \star =id,&\label{co4}
		\end{eqnarray}
                where $\Sigma:A\otimes A\to A\otimes A$ is the flip
                $\Sigma(a_1\otimes a_2)=a_2\otimes a_1$;
	\item [(HS$_2$)]
                the mapping $\Delta:A\to A\otimes A$ is
		positive, i.e. it maps the cone of positive elements
		of $A$ into the cone of positive elements of
		$A\otimes A$;
	\item[(HS$_3$)]
		the following identities hold
		\begin{eqnarray}
			&
			(a\cdot b)^\star=a^\star\cdot b^\star ,\qquad
                        \Delta\circ *=(*\otimes *)\circ \Delta,
			&\label{hy:1}\\[2mm]
			&
			\epsilon(a\cdot b)=\epsilon(a)\epsilon(b),
			\qquad
                        \Delta(1)=1\otimes 1,
			&\label{hy:2}\\[2mm]
			&
			\star \circ *=*\circ \star.
			&\label{hy:3}
		\end{eqnarray}
\end{itemize}

{\bf 1.2.} Let $(A,\Delta,\epsilon,\star)$ be a hypergroup structure on
a $C^*$-algebra $A$. By $A'$ we denote the set of all continuous
linear functionals on the $C^*$-algebra $A$. For $\xi,\eta\in A' $
 we define a product $\cdot $ and an involution $^+$ by
\begin{equation}\label{fun_alg}
	\begin{array}{c}
		(\xi\cdot \eta)(a)=(\xi\otimes \eta)\delta(a),\\[2mm]
		\xi^+(a)=\overline {\xi(a^\star )},
	\end{array}
\end{equation}
$a\in A$, with the norm given by
\begin{equation}\label{fun_norm}
	\|\xi\|=\sup_{\|a\|=1}|\xi(a)|.
\end{equation}
Then $A'$ is a unital Banach $\ast$-algebra.

A state $\mu\in A' $ is
called a {\em Haar measure} (with respect to the
hypergroup structure) if
	\begin{equation}\label{inv_int}
                (\mu\otimes id)\circ \Delta(a)=(id\otimes
                \mu)\circ \Delta(a)=\mu(a)1
	\end{equation}
	for all $a\in A$.

Let $(A,\Delta,\epsilon,\star)$ be a hypergroup structure on
	the $C^*$-algebra $A$. An element $a\in A$ is called
	{\em positive definite} if
	\begin{equation}\label{pos_def}
		\xi\cdot \xi^+(a)\geq 0
	\end{equation}
	for all $\xi\in A'$.

{\bf 1.3. Theorem} (\cite{ChaVai}).
       {\em Let $(A,\Delta,\epsilon,\star)$ be a hypergroup structure on
	a $C^*$-algebra $A$. Suppose that the linear space spanned
	by positive definite elements is dense in $A$. Then there
        exists a Haar measure $\mu$, it is unique, and
        $\mu^+=\mu$.}

{\bf 1.4.} Suppose that $(A,\Delta,\epsilon,\star )$ is a hypergroup
	structure on a $C^*$-alge\-bra $(A,\cdot ,1,*)$ and $\mu$ is
	the corresponding Haar measure. We call
        ${\cal A}=(A,\Delta,\epsilon,\star ,\sigma_t,\mu)$ a
	{\em compact quantum hypergroup} if
	\begin{itemize}
	\item [($QH_1$)] the mapping $\Delta$ is completely positive
		and the linear span of positive
		definite elements is dense in $A$;
	\item [($QH_2$)]
		$\sigma_t$, $t\in {\Bbb R}$, is a continuous
		one-parameter group of automorphisms of $A$ such that
		\begin{itemize}
		\item [(a)]
		        there exist dense subalgebras $A_0 \subset A$ and
          		$\tilde{A_0} \subset A\otimes A$ such that the one-parameter
                        groups $\sigma_t$ and $\sigma_t\otimes \ \mbox{id}$,
                        $\ \mbox{id}\otimes \sigma_t$  can
			be extended to complex one-parameter
                        groups $\sigma_z$ and $\sigma_z\otimes \ \mbox{id}$,
                        $\ \mbox{id}\otimes \sigma_z$, $z\in{\Bbb C}$, of
	        	automorphisms of the algebras $A_0$ and $\tilde{A_0}$
	        	respectively;
		\item [(b)]
			$A_0$ is invariant with respect to $*$ and
                        $\star $, and $\Delta(A_0)\subset \tilde{A_0}$;
		\item [(c)]
			the following relations hold on $A_0$ for all
			$z\in {\Bbb C}$:
			\begin{eqnarray}
				&\displaystyle
                                \Delta\circ \sigma_z=(\sigma_z\otimes
                                \sigma_z)\circ \Delta,
				&\label{s_d}\\
				&\displaystyle
                                \mu(\sigma_z(a))=\mu(a);\label{s_inv_int}
				&
			\end{eqnarray}
		\item [(d)]
			there exists $z_0\in {\Bbb C}$ such that the
                        Haar measure $\mu$ satisfies the following
			strong invariance condition for all $a,b\in
			A_0$:
			\begin{equation}\label{s_inv}
                                (id \otimes \mu)\bigl[\bigl(( *\circ
				\sigma_{z_0}\circ \star \otimes
                                id)\circ \Delta(a)\bigr)\cdot(1\otimes
				b)\bigr] =(id\otimes
                                \mu)\bigr((1\otimes
                                a)\cdot\Delta(b)\bigr);
			\end{equation}
		\end{itemize}
	\item [($QH_3$)]
                the Haar measure $\mu$ is faithful on $A_0$.
	\end{itemize}

        It will be convenient to denote
\begin{equation}\label{k}
	\kappa=*\circ \sigma_{z_0}\circ \star
\end{equation}
and call it an {\em antipode}. Then for all $a, b \in A_0$,
	\begin{equation}\label{k_prop}
		\begin{array}{c}
		\kappa(ab)=\kappa(b)\kappa(a),\qquad \delta\circ
		\kappa=\Pi\circ (\kappa\otimes \kappa)\circ \delta,
		\\ [2mm]
		\nu\circ \kappa=\nu,\qquad \kappa(1)=1,\qquad
		\epsilon\circ \kappa=\epsilon.
		\end{array}
        \end{equation}
Note that $\kappa$ is invertible with $\kappa^{-1}=\star \circ
\sigma_{-z_0}\circ *$.

With such a notation, relation (\ref {s_inv}) becomes
\begin{equation}\label{s_inv_int'}
        (id \otimes \mu)\bigl(( \kappa
	\otimes
        id)\circ \Delta(a)\cdot(1\otimes
	b)\bigr) =(id\otimes
        \mu)\bigr((1\otimes
        a)\cdot\Delta(b)\bigr).
\end{equation}

\par {\bf 1.5.} Compact quantum groups are quantum hypergroups.
In detail, let ${\cal A} = (A,\Delta, \epsilon, \kappa)$ be a compact matrix pseudogroup
with $A_0$ being the involutive subalgebra generated by
matrix elements of the fundamental corepresentation
\cite{Wor}. We will use the following notations
\begin{equation}\label{conv}
        \xi\cdot a=(id\otimes \xi)\circ \Delta(a),\qquad
        a\cdot \xi=(\xi\otimes id)\circ \Delta(a),\qquad
	\xi\cdot \eta=(\xi\otimes \eta)\circ \Delta
\end{equation}
for $\xi,\eta\in A'$ and $a\in A$.
It readily follows from (\ref {conv}) that
\begin{equation}\label{conv_ass}
	\xi\,.\,(\eta\,.\,a)=(\xi\cdot \eta)\,.\,a,\qquad
	(a\,.\,\eta)\,.\,\xi=a\,.\,(\eta\cdot \xi).
\end{equation}

	Let $U^\alpha=(u^\alpha_{ij})_{i,j=1}^{d_\alpha}$ be an
irreducible unitary corepresentation of $A$. Then there exists a
unique, up to a positive constant, positive definite matrix
$M^\alpha=(m^\alpha_{ij})_{i,j=1}^{d_\alpha}$ such that
\begin{equation}\label{def:F}
	M^\alpha\cdot U^\alpha=S^2(U^\alpha)\cdot M^\alpha,
\end{equation}
where $\cdot $ here denotes the usual matrix multiplication \cite{Wor}.

	For each $z\in {\Bbb C}$, we denote by $m_{ij}^{\alpha (z)}$
the matrix elements of the matrix $(M^\alpha)^z$. It is known that
there exists a one-parameter family of homomorphisms $f_z:A_0\to{\Bbb
C}$, $z\in {\Bbb C}$, where, as before, $A_0$ denotes the
$*$-subalgebra generated by matrix elements of the fundamental
corepresentation. These homomorphisms are defined
by
\begin{equation}\label{def:f}
	f_z(u^\alpha_{ij})=m^{\alpha(z)}_{ij}
\end{equation}
and possess the following properties \cite{Wor}:
\begin{itemize}
	\item [$(F_1)$:]
		$f_z(1)=1$ for all $z\in {\Bbb C}$;
	\item [$(F_2)$:]
		$f_z\cdot f_{z'}=f_{z+z'}$ and $f_0=\epsilon$;
	\item [$(F_3)$:]
		$f_z(\kappa(a))=f_{-z}(a)$;
	\item [$(F_4)$:]
		$f_z(a^*)=\overline {f_{-\overline z}(a)}$;
	\item [$(F_5)$:]
		$\kappa^2(a)=f_{-1}\cdot a\cdot f_1$;
	\item [$(F_6)$:]
                $\mu(a\cdot b)=\mu(b\cdot (f_1\cdot a\cdot f_1))$,
                where $\mu$ is the Haar measure on the compact
		quantum group $A$.
\end{itemize}

        It was shown in \cite{ChaVai} that
$(A, \Delta, \epsilon, \star, \sigma_t)$ is
quantum hypergroup,  where the mapping $\star :A_0\to A_0$ is defined by
\begin{equation}\label{def:star}
        a^\star =f_{-1/2}\cdot S(a)^*\cdot f_{1/2},
\end{equation}
the action of the group $\sigma_t$ is defined by
\begin{equation}\label{Wor:sigma}
\sigma_t (a) = f_{it}\cdot a\cdot f_{-it}
\end{equation}
and $z_0 = -\frac 1 2 i$.

{\bf 1.6.} Here we briefly discuss the case of usual compact hypergroups,
or equivalently, normal hypercomplex systems with basis unity
(\cite{BK}). Let ${\cal A}$ be a commutative
compact quantum hypergroup. Let $Q$ denote
the spectrum of the commutative $C^*$-algebra $A$.
Each element $\xi\in Q$ defines a linear
operator on $A$ given by
\begin{equation}\label{def:gto}
        R_\xi=(id\otimes \xi)\circ \Delta.
\end{equation}
The operators $R_\xi$, $\xi\in Q$, constitute a family of
generalized translation
operators on $A\cong C(Q)$.
For $\xi\in Q$ and $a\in A_0$, we define
\begin{equation}\label{char*}
        \xi^\dagger(a)=\overline {\xi(a^\dagger)},
\end{equation}
where
\begin{equation}\label{a*}
		a^\dagger=\kappa(a)^* .
\end{equation}
 From the definition of $a^\dagger$, it immediately follows that
$\xi^\dagger$ is a homomorphism $A_0\to {\Bbb C}$ and, hence,
continuous. Being extended by continuity to $A$, it becomes a point
in $Q$. In \cite{ChaVai} it was established that $Q$ is a basis of a normal
hypercomplex system $L_1(Q,\mu)$ with a basis unit $\epsilon$.

\par {\bf 1.7. Preliminary on twisting of Kac algebras.}
We recall the construction of a twisting of Kac algebras following
\cite{EnVai}, \cite{Vai:twist}.
For our purposes, it suffices to deal only with the case of finite
dimensional Kac algebras.

Let ${\cal A} = (A,
\Delta, \varepsilon, \kappa)$ be a finite dimensional Kac algebra.  A {\it
2-cocycle} of $A$ is a unitary $\Omega$ in $A\otimes A$ such that
$$(\Omega\otimes 1)(\Delta\otimes
\ \mbox{id})(\Omega)=(1\otimes\Omega)(\ \mbox{id}\otimes
\Delta)(\Omega).$$
A {\it 2-pseudo-cocycle} of $A$ is a unitary $\Omega$ in
$A\otimes A$ such that
$$
\partial_2 \Omega = (\ \mbox{id}\otimes\Delta)(\Omega^*)
(1\otimes\Omega^*)(\Omega\otimes 1)(\Delta\otimes \ \mbox{id})(\Omega)\
$$
belongs to $(\Delta\otimes \ \mbox{id})\Delta (A)'$ (\cite{EnVai}, 2.3).
Let us put, for all $x$ in $A$ :
$$\Delta_\Omega (x)=\Omega\Delta (x)\Omega^* .$$
Then $\Delta_\Omega$ is coassociative iff $\Omega$
is a 2-pseudo-cocycle of $A$; and we shall say that
$\Delta_\Omega$ is {\it twisted} from $\Delta$ by $\Omega$.

\par For a unitary $u\in A$, let us put
$\Omega^u = (u^*\otimes u^*)\Omega\Delta(u)$. Then
$\Omega^u$ is a 2(-pseudo)-cocycle iff $\Omega$ is.
Let $\Sigma(x\otimes y)=y\otimes x$ be the flip in $A\otimes A$.
We shall say that a 2(-pseudo)-cocycle
$\Omega$ is {\sl pseudo-co-involutive} (resp., {\sl co-involutive}) if
$\Omega^u(\Sigma(\kappa\otimes\kappa)(\Omega^*))^*\in \Delta(A)'$
(resp.,$\Sigma(\kappa\otimes\kappa)(\Omega^*)=
\Omega^u)$. We shall say that $\Omega$ is {\sl strongly
co-involutive}, if $\Sigma(\kappa\otimes\kappa)(\Omega^*)=\Omega$.
It should be noted that each {\sl counital} 2-cocycle $\Omega$ of
a finite dimensional Kac algebra (i. e. such that
$(\varepsilon\otimes\ \mbox{id})\Omega = (\ \mbox{id}\otimes\varepsilon)\Omega = 1$)
is coinvolutive with $u = m(\ \mbox{id}\otimes\kappa)\Omega = \kappa (u)$ and
$\varepsilon(u) = 1$ (here $m$ denotes the multiplication
in the algebra $A$).
If $\Omega$ is a (pseudo)-co-involutive 2(-pseudo)-cocycle of a Kac algebra $A$,
then the coalgebra $(A, \Delta_\Omega)$ possesses a coinvolution
of the form $\kappa_\Omega (x) = u\kappa(x)u^*$. If
$\Omega$ is a strongly coinvolutive 2-cocycle of a Kac algebra $A$,
then $\kappa$ is a coinvolution of the coalgebra $(A, \Delta_\Omega)$.

\par Let $G$ be a finite group and $H$ be its Abelian subgroup.
Let ${\Bbb C}(G)$ be the $C^*$-algebra generated by the left regular
representation $\lambda$ of $G$. Then ${\Bbb C}(G)$ has a standart structure of
cocommutative Kac algebra $({\Bbb C}(G), \Delta, \kappa, \mu)$, where
$\Delta(\lambda(g)) = \lambda(g)\otimes\lambda(g)$, $\kappa(\lambda(g)) =
\lambda(g^{-1})$,
$\mu(\lambda(f)) = f(e)$,  $g\in G$ and
$f$ is a continuous function on $G$. Denote by $\hat H$ the dual group of $H$.
Then there exists a family of projections $P_{\hat h}$
generating an Abelian subalgebra ${\Bbb C}(H)$ in ${\Bbb C} (G)$ such that
$$P_{\hat h}P_{\hat g}=\delta_{\hat h, \hat g}P_{\hat h},\
\lambda (h)=\sum_{\hat h\in \hat H}<\hat h, h>P_{\hat h},\
P_{\hat h}={1\over {||H||}}\sum_{h\in H}<\hat h, h>\lambda (h),
$$
$$
\Delta (P_{\hat h})=\sum_{\hat g\in \hat H}P_{\hat g}\otimes
P_{{\hat g}^{-1}\hat h},\ \kappa(P_{\hat h})=P_{{\hat h}^{-1}},
$$
where $g,h\in H,\ \hat g, \hat h\in \hat H,\ ||H||$ is the order of
$H$.
Using these idempotents, one can write the following formulae for
$\Omega$ and $u$:
$$\Omega=\sum_{x,y\in\hat H}\omega(x,y)(P_x\otimes P_y),\
u=m(i\otimes\kappa)\Omega=\sum_{x\in\hat H}\omega(x,x^{-1})P_x,$$
where  $\omega:\hat H \times \hat H \to{\Bbb T}$ is
a 2-(pseudo)-cocycle.
If $\omega$ is a 2-cocycle, then so is $\Omega$, but if
$\omega$ is only a 2-pseudo-cocycle on $\hat H$,
one should verify that $\Omega$ is a 2-pseudo-cocycle on
$({\Bbb C}(G),\ \Delta,\ \kappa)$. If it is, and it
is at least pseudo-coinvolutive, then we have already a new finite
dimensional coinvolutive coalgebra. If it has a counit, then it is
a Kac algebra. For this, one should choose $\omega$ to be
counital, $\omega(e,x)=\omega(x,e)=1\ \forall x\in H$, which gives the
counitality of $\Omega$.

\subsection*{2. Conditional expectations on quantum hypergroups}
\setcounter{subsection}{2}
\setcounter{equation}{0}

Let $(A,\Delta,\epsilon,\star ,\sigma_t,\mu)$ be a compact
quantum hypergroup and $\kappa$ be the antipod defined by (\ref{k}).
Let $B$ be a unital $C^*$-subalgebra of $A$ and $P: A \to B$ be
a corresponding $\mu$-invariant conditional expectation \cite{cond_exp}
(i. e. $\mu\circ P = \mu$). The following Theorem 2.1
states that, under some conditions, we can define, on the algebra $B$,
a new comultiplication $\tilde\Delta$ such that
$(B, \tilde\Delta, \epsilon, \star,\sigma_t, \mu)$
is a compact quantum hypergroup.

\par {\bf 2.1. Theorem.} {\em Let ${\cal A}$ be a compact quantum hypergroup,
           $B$ an unital $C^*$-subalgebra of a $C^*$-algebra $A$ and
           $P: A \to B$ a $\mu$-invariant conditional expectation.
           Then, let us put for all $x \in B$,
\begin{equation}\label{_Delta}
             \tilde\Delta (x) = (P\otimes P) \Delta (x).
\end{equation}
          Suppose that the following conditions hold:
         \begin{itemize}
         \item [(1)] $(P\otimes P)\Delta (x) = (P\otimes P)\Delta (P(x))$
                     (or, equivalently, $\ker P$ is a coideal in $A$),
         \item [(2)] $P\circ\star = \star\circ P$,
         \item [(3)] the dense subalgebras $A_0 \subset A$ and
                     $\tilde{A_0}\subset A\otimes A$ are invariant under
                     $P$ and $P\otimes P$, respectively, and $P\circ\sigma_z =
                     \sigma_z\circ P$ for all $z\in {\Bbb C}$,
         \item [(4)] the mapping $\epsilon$ is a counit of the coalgebra
                     $(B, \tilde\Delta )$, i. e., the relation
                     $$
                     (\epsilon\otimes id)\circ \Delta=(id\otimes \epsilon)\circ \Delta=id
                     $$
                     holds.
         \end{itemize}
         Then $(B, \tilde\Delta, \varepsilon,\star ,\sigma_t,\mu)$ is a compact quantum
         hypergroup.
         }
\par {\bf Proof.} First we prove that $(B, \tilde\Delta, \epsilon,\star )$ is
a hypergroup structure. Using the first condition of the theorem, we
have for all $x\in B$
\begin{eqnarray*}
  (\tilde\Delta\otimes \ \mbox{id})\tilde\Delta (x)
  &=& (P\otimes P\otimes \ \mbox{id})(\Delta\otimes \ \mbox{id})(P\otimes P)\Delta(x)\\
  &=& (\ \mbox{id}\otimes \ \mbox{id}\otimes P)((P\otimes P) \circ \Delta \circ P \otimes \ \mbox{id})\Delta(x) \\
  &=& (P\otimes P\otimes P)(\Delta\otimes \ \mbox{id})\Delta(x)\\
  &=& (P\otimes P\otimes P)(\ \mbox{id}\otimes \Delta)\Delta(x)
  = (\ \mbox{id} \otimes \tilde\Delta)\tilde\Delta (x),
\end{eqnarray*}
i.e., the mapping $\tilde\Delta : B \to B\otimes B$ is coassociative.
Relation (\ref{co2}) follows from Condition 4 and relations
(\ref{co3}), (\ref{co4}) easily follow from the second condition of the theorem.
The mapping $\tilde\Delta $ is completely positive
since $P$ is completely positive and $\Delta$ is a homomorphism.
Finally, identities (\ref{hy:1})-(\ref{hy:3}) are obvious.

\par Let $A',\ B' $ be the spaces of continuous linear functionals
on $A,\ B$, respectively, equipped with the multiplications
$$
\left<\xi\cdot\eta, a\right> = \left<\xi\otimes\eta, \Delta (a)\right>\qquad (\xi,\eta \in A'; a \in A)
$$
and
$$
\left<\xi'\cdot\eta', b\right> = \left<\xi'\otimes\eta', \tilde\Delta (b)\right> \qquad (\xi',\eta' \in B'; b \in B).
$$
In order to show that the linear span  of positive definite elements
is dense in $B$, let us define a mapping $P': B' \to A'$
as follows:
\begin{equation}\label{P'}
\left< P'\xi,\ a\right> = \left< \xi,\ Pa\right>,
\end{equation}
 where $a\in A, \xi \in B'$.
If $a\in A$ is a positive definite element then $Pa\in B$ is a positive
definite element (with respect to $\tilde\Delta$). Indeed, for all
$\xi\in B'$ we have
\begin{eqnarray*}
\left<\xi\cdot\xi^+, P a\right>
&=& \left<\xi\otimes\overline{\xi}, (id\otimes \star)(P\otimes P)\Delta (P a)\right>\\
&=& \left<\xi\otimes\overline{\xi}, (P\otimes P)(id\otimes \star)\Delta (a)\right>\\
&=& \left<P' \xi\otimes\overline{P' \xi}, (id\otimes \star)\Delta (a)\right> \geq 0,
\end{eqnarray*}
because $\left<\eta\cdot \eta^+, a\right> =
\left<\eta\otimes\overline{\eta}, (id\otimes \star)\Delta (a)\right> \geq 0$
for all $\eta\in A'$. So the axiom ($QH_1$) holds for the hypergroup structure
$(B, \tilde\Delta, \epsilon,\star )$. Since the linear span of
$\{P a \in B\ |\ a \ \mbox{is positive definite in}\  A\}$ is dense in $B$,
we get the result.

\par Define $B_0 = P A_0$ and $\tilde{B_0} = P\otimes P (\tilde{A_0})$.
Conditions 2 and 3 of the theorem imply that $B_0$ and $\tilde{B_0}$ are
dense subalgebras of $B$ and $B\otimes B$, respectively, such that
axioms ($QH_2$)(a) - ($QH_2$)(c) hold.

\par Let us check the remaining axioms. The faithful state $\mu$ is a Haar measure
of the hypergroup structure $(B,\tilde\Delta, \epsilon, \star) $,  since
$P$ is $\mu$-invariant.
Let us show that $\mu$ is strongly invariant with respect to $\tilde\Delta$,
i.e. relation (\ref{s_inv_int'}) holds for the hypergroup structure
$(B,\tilde\Delta, \epsilon, \star)$. It follows from Conditions 2 and 3
of the theorem, that $P$ commutes with the antipode $\kappa $ defined by (\ref{k}).
Using this fact, the fact that a Haar measure in a quantum hypergroup ${\cal A}$
is strongly invariant, and since $P$ is $\mu$-invariant, we have for all $a, b \in B$
that
\begin{eqnarray*}
& &(id\otimes\mu)[((\kappa\otimes id)\circ\tilde\Delta (a))\cdot(1\otimes b)]\\
&=& (id\otimes\mu)[((P\otimes P)\circ(\kappa\otimes id)\circ\Delta (a))\cdot(1\otimes b)] \\
&=& (id\otimes\mu)\circ (P\otimes P)[((\kappa\otimes id)\circ\Delta (a))\cdot(1\otimes b)]\\
&=& (P\otimes\mu)[((\kappa\otimes id)\circ\Delta (a))\cdot(1\otimes b)] \\
&=& (P\otimes\mu)[(1\otimes a)\cdot\Delta (b))]
= (id\otimes\mu)[(P\otimes P)((1\otimes a)\cdot\Delta (b)] \\
&=& (id\otimes\mu)[(1\otimes a)\cdot\tilde\Delta (b)].
\end{eqnarray*}\hfill$\Box$

\par {\bf 2.2. Remark.} We can simplify the hypothesis of Theorem 2.1
in case when ${\cal A}$ is a compact quantum group. Indeed,
let ${\cal A} = (A, \Delta, \epsilon, \kappa)$ be a compact matrix
pseudogroup with $A_0$ being the $\ast$-subalgebra
generated by matrix elements of the fundamental corepresentation.
Let $\mu$ be its Haar measure and $P: A\to B$ be a $\mu$-invariant
conditional expectation that maps to a unital $C^*$-subalgebra $B$ of $A$.
Let us define a new comultiplication $\tilde\Delta$ on $B$ in accordance with
(\ref{_Delta}). Then the statement of  Theorem 2.1 remains true if
we replace the hypothesis 2 and 3 with the following condition:
\begin{itemize}
   \item [(A)] $A_0$ is invariant with respect to $P$ and
              $\kappa \circ P = P\circ\kappa$.
\end{itemize}

\par {\bf Proof.} In this case we have that $\tilde{A_0}\subset A_0\otimes A_0$
and we only need to prove that $P$ commutes with $\star$ and $\sigma_t$.
Using property $(F_5)$ of $f_z$ we have for all $a\in A_0$,
$$
P(f_{-1}\cdot a\cdot f_{1}) = P(\kappa^2 (a)) = \kappa^2 (P(a)) = f_{-1}\cdot Pa\cdot f_1,
$$
whence $P(f_{-n}\cdot a\cdot f_n) = f_{-n}\cdot Pa\cdot f_n$ for all $n\in {\Bbb N}$.
Since $f_z (a)$ is an entire function of exponential growth on the right half-plane
for all $a\in A_0$, the vector-valued function
$$
f_{-z}\cdot u^\alpha_{kl}\cdot f_z = \sum_{r, s} f_{z}(\kappa(u^\alpha_{kr})) u^\alpha_{rs} f_z (u^\alpha_{sl})
$$
is entire of exponential growth on the right half-plane in the weak sence,
and we obtain (cf. \cite{Wor}, Lemma 5.5) that $P(f_{-z}\cdot a\cdot f_z) =
(f_{-z}\cdot Pa\cdot f_z)$ for all $z\in {\Bbb C}$ and $a\in A_0  $,  whence
it follows from (\ref{Wor:sigma}) that $P$ commutes with $\sigma_z$, $z\in {\Bbb C}$.
Also, for $a\in A_0$ by using (\ref{def:star}), we have
$P(a^\star) = P(f_{-1/2}\cdot \kappa(a)^*\cdot f_{1/2}) =
f_{-1/2}\cdot P((\kappa(a))^*)\cdot f_{1/2} = f_{-1/2}\cdot (\kappa(Pa))^*\cdot f_{1/2} =
(Pa)^\star$.
\hfill$\Box$

\par In order to clarify Condition 4 of  Theorem 2.1,
we establish sufficient
conditions for the mapping $\epsilon $ be a counit of the
coalgebra $(B,\tilde\Delta)$. The following Remark 2.3
is quite transparent, but its condition does not hold
in some examples. Therefore, we establish Proposition 2.4
with less restrictive conditions on $\epsilon$
in the case where ${\cal A}$ is a finite dimensional Kac algebra
(finite dimensional quantum groups are, in fact,
Kac algebras \cite{Wor}).

\par {\bf 2.3. Remark.}  Let a $\mu$-invariant
      conditional expectation $P: A\to B$ on a
      compact quantum hypergroup ${\cal A}$ satisfy Condition 1 of Theorem 2.1.
      If the mapping $\epsilon : A \to {\Bbb C}$ satisfies the relation
      $\epsilon \circ P = \epsilon$, then the mapping
      $\epsilon : B\to {\Bbb C}$ is a counit of the
      coalgebra $(B, \tilde\Delta)$, i. e. (\ref{co2}) holds.


\par {\bf 2.4. Proposition.} {\sl Let a $\mu$-invariant
      conditional expectation $P: A\to B$ on a finite dimensional
      Kac algebra ${\cal A}$ satisfy  Condition 1 of  Theorem 2.1.
      Denote by $p_{\epsilon}$ the one-dimensional central projection
      in $A$ such that the relation $a p_{\epsilon} = \epsilon (a) p_{\epsilon}$
      holds for all $a\in A$.
      Let the following conditions hold for all $b\in B$
      \begin{eqnarray} \label{p_epsilon1}
      (P\otimes P) [(p_{\epsilon}\otimes 1)\cdot\Delta(b)]
          &=&(P\otimes P) [(P(p_{\epsilon})\otimes 1)\cdot\Delta(b)]\\
      \label{p_epsilon2}
      (P\otimes P) [(1 \otimes p_{\epsilon})\cdot\Delta(b)]
          &=&(P\otimes P) [(1 \otimes P(p_{\epsilon}))\cdot\Delta(b)].
      \end{eqnarray}
      Then the mapping
      $\epsilon : B\to {\Bbb C}$ is a counit of the finite dimensional
      coalgebra $(B, \tilde\Delta)$.

}

\par {\bf Proof.} First we prove that relations
(\ref{co2}) hold if and only if, for all $b\in B$, we have
\begin{eqnarray}\label{p_eps1}
(P(p_{\epsilon})\otimes 1)\tilde\Delta (b) &=& P(p_{\epsilon})\otimes b,\\
\label{p_eps2}
(1 \otimes P(p_{\epsilon}))\tilde\Delta (b) &=& b\otimes P(p_{\epsilon}).
\end{eqnarray}
Indeed, from (\ref{p_eps1}) with the usual tensor notation $\tilde\Delta (b) = b_{(1)}\otimes b_{(2)}$,
we have
\begin{eqnarray*}
P(p_{\epsilon})\otimes b &=& (P(p_{\epsilon})\otimes 1)\tilde\Delta (b)
= P(p_{\epsilon} b_{(1)})\otimes b_{(2)} \\
&=& \epsilon (b_{(1)})P(p_{\epsilon})\otimes b_{(2)}
= P(p_{\epsilon})\otimes \epsilon (b_{(1)})b_{(2)}\\
&=& P(p_{\epsilon})\otimes (\epsilon\otimes id)\tilde\Delta (b).
\end{eqnarray*}
Similarly, we have
$b \otimes P(p_{\epsilon}) = (id\otimes \epsilon)\tilde\Delta (b)\otimes P(p_{\epsilon}).$
The converse statement follows by similar arguments.

\par Now by using (\ref{p_epsilon1}) we get (\ref{p_eps1}):
\begin{eqnarray*}
(P(p_{\epsilon})\otimes 1)\tilde\Delta (b)
&=& (P\otimes P) [(P(p_{\epsilon})\otimes 1)\cdot\Delta(b)]\\
&=& (P\otimes P) [(p_{\epsilon}\otimes 1)\cdot\Delta(b)]
= (P\otimes P) [(p_{\epsilon}\otimes b]\\
&=& P(p_{\epsilon})\otimes b.
\end{eqnarray*}
Similarly, (\ref{p_eps2}) follows from (\ref{p_epsilon2}).\hfill$\Box$

\par {\bf 2.5. Lemma.} {\sl Let ${\cal A}$ be a compact quantum hypergroup
and all hypothesis of  Theorem 2.1 hold for a conditional expectation
$P: A\to B$. Let $P': B'\to A'$ be defined by (\ref{P'}). Then
$P'$ is a $\ast$-homomorphism.
}

\par {\bf Proof.} For any $\eta, \eta_1, \eta_2 \in B'$ and $a\in A$, we have
\begin{eqnarray*}
\left< P'(\eta_1\cdot\eta_2), a\right> &=& \left< \eta_1\cdot\eta_2, Pa\right>
= \left< \eta_1\otimes\eta_2, \tilde\Delta (Pa)\right>\\
&=& \left< \eta_1\otimes\eta_2, (P\otimes P) \Delta (Pa)\right>
= \left< \eta_1\otimes\eta_2, (P\otimes P) \Delta (a)\right> \\
&=& \left< P'(\eta_1)\otimes P'(\eta_2), \Delta (a)\right>
= \left< P'(\eta_1)\cdot P'(\eta_2), a\right>;\\
\left< P'(\eta^+), a\right> &=& \overline{\left< \eta, P(a)^\star\right>}
=\overline{\left< \eta, P(a^\star)\right>} = \left< (P'\eta)^+, a\right>.
\end{eqnarray*}
\hfill$\Box$

\par {\bf 2.6.} In what follows we discuss conditional expectations
on usual hypergroups and a construction of
orbital morphisms \cite{Jew}. For simplicity, we consider only
the case of compact hypergroups. Let $Q$ be a compact
DJS-hypergroup with involution $Q\ni p\mapsto p^\dagger \in Q$,
comultiplication $\Delta : C(Q) \to C(Q)\otimes C(Q)$,
neutral element (counit) $e$ and Haar
measure $\mu$, and let $Y$ be a compact Hausdorff space.
Let $\phi$ be an open continuous mapping from $Q$ onto $Y$
(orbital mapping).
The closed sets $\phi^{-1} (y)$, $y\in Y$,
are called {\sl $\phi$-orbits}.
Let $B$ be a $C^*$-subalgebra of
$C(Q)$ consisting of functions constant on $\phi$-orbits. Obviously,
a mapping $\phi : C(Y) \to B$ defined by $(\phi f)(x) = f(\phi (x))$
is an isomorphism of the $C^*$-algebras.
Denote by $\phi_\ast : M(Q) \to M(Y)$
the corresponding mapping of Radon measures,
$\left< \phi_\ast (\mu), f\right> =
\left< \mu, f\circ \phi \right>$, $\mu\in M(Q),\ f\in C(Y)$.
The measure $\mu\in M(Q)$ is called {\sl $\phi$-consistent},
if $\phi_\ast (\mu\ast\nu) = 0 = \phi_\ast (\nu\ast\mu)$ whenever
$\phi_\ast (\nu) = 0$. The following proposition clarifies
the concept of $\phi$-consistency.

\par {\bf 2.7. Proposition.} {\sl Let $\phi: Q\to Y$ be an orbital mapping and
denote $\tilde e = \phi (e)$. Suppose that the following conditions
are satisfied:
\begin{itemize}
  \item[$(a)$] $\phi^{-1}(\tilde e) = \{e\}$,
  \item[$(b)$] if $A$ is $\phi$-orbit then so is $A^\dagger$,
  \item[$(c)$] for any $y\in Y$ there exists a probability measure
             $q_y \in M(Q)$ such that $\ \mbox{supp}\ q_y \subset \phi^{-1}(y)$,
  \item[$(d)$] each measure $q_y$ is $\phi$-consistent.
\end{itemize}
Define the linear mapping $P: C(Q) \to B$ as follows:
\begin{equation}\label{P:orb}
(P f)(x) = \left< q_{\phi(x)}, f\right>.
\end{equation}
Then $P$ is a conditional expectation and satisfies all hypothesis
of Theorem 2.1. Conversely, if conditions $(a)$--$(c)$ hold and
the linear mapping $P$ defined by (\ref{P:orb}) is a conditional expectation
from $C(Q)$ to $B$  satisfying all
hypothesis of Theorem 2.1, then each $q_y$ is $\phi$-consistent.}

\par {\bf Proof.} By virtue of $(b)$ we can define an involutive
homeomorphism $\dagger: Y\to Y$ as follows: if $y=\phi(x)$, then
$y^\dagger = \phi (x^\dagger)$.
Theorem 13.5A in \cite{Jew} states that
there exists a unique convolution $\ast$ in $M(Y)$ such that $Y$
is a hypergroup and $\phi$ is an orbital morphism, i. e.,
\begin{itemize}
       \item [(i)] $\delta_y \ast \delta_z = \phi_* (q_y\ast q_z)$
                   for any $y, z \in Y$,
       \item [(ii)] $q_{y^\dagger} = (q_y)^\dagger$,
       \item [(iii)] $\ \mbox{supp}\ q_y = \phi^{-1}(y)$ and
                    \begin{equation}\label{recomp}
                    m = \int_Q q_{\phi(x)} m(dx),
                    \end{equation}
                    where $m$ is a Haar measure on the hypergroup $Q$.
       \end{itemize}
Thus we can define a mapping $\phi^*: M(Y)\to M(Q)$ by setting
$\phi^* (\delta_y) = q_y$ and
$$
\phi^* (\nu) = \int_Y q_z \  \nu(dz),
$$
for $y\in Y$ and $\nu\in M(Y)$. In virtue of Lemma 13.6A in \cite{Jew},
the orbital morphism $\phi$ is consistent, i. e.
the mapping $\phi^*$ is a $\ast$-homomorphism.
It also follows from the proof of Theorem 13.5A cited above that
the mapping $Y\ni y\mapsto q_y \in M(q)$ is continuous in the weak topology.
Thus $P$ is well-defined and is, indeed, a conditional expectation.
Let $f\in \ker P$. Then $\left<q_z, f\right> = 0$ for all $z\in Y$, and
for $x_1, x_2 \in Q$ we have
\begin{eqnarray*}
((P\otimes P) \Delta f)(x_1, x_2) &=&
\left<q_{\phi(x_1)}\otimes q_{\phi(x_2)}, \Delta f\right>\\
&=& \left<q_{\phi(x_1)}\ast q_{\phi(x_2)}, f\right>
= \left<\phi^*(\delta_{\phi(x_1)})\ast \phi^*(\delta_{\phi(x_2)}), f\right> \\
&=& \left<\phi^*(\delta_{\phi(x_1)}\ast \delta_{\phi(x_2)}), f\right>
= \int_Y \left<q_z, f\right>\ (\delta_{\phi(x_1)}\ast \delta_{\phi(x_2)})(dz) \\
&=& \int_Y (Pf)(\phi(x)) \ (\delta_{\phi(x_1)}\ast \delta_{\phi(x_2)})(dz)=0,
\end{eqnarray*}
where $x\in\phi^{-1}(z)$. Hence, $\ker P$ is a coideal.
In virtue of (iii), $P$ is $m$-invariant and the fourth condition of Theorem 2.1
follows from $(a)$. At last, the equality $P\circ\dagger = \dagger\circ P$
follows from $(b)$ (note that, for usual hypergroups,
$f^\star (x) = \overline{f(x^\dagger)}$ and $\sigma_{z} = id$).

\par Let us prove the converse statement.
Define the convolution in $M(Y)$ in accordance with (i).
Since $B$ is isomorfic to $C(Y)$, it follows from Theorem 2.1
that $Y$ is a DJS-hypergroup.
To prove the result, we need to show that $\phi$ is
a consistent orbital morphism. Then the result follows from
Theorem 13.6B in \cite{Jew}. Indeed, (ii) follows from $(b)$
and equality (\ref{recomp}) follows from the fact that $P$ is
$m$-invariant. Let us show that $\ \mbox{supp}\ q_y = \phi^{-1} (y)$. Suppose that
$x_0\in \phi^{-1} (y)$ and $x_0 \not\in \ \mbox{supp}\ q_y$. Then
there exists an open neighborhood $O_{x_0}$ of $x_0$ such that
$O_{x_0}\cap\ \mbox{supp}\ q_y = \varnothing$. Let $f\in C(Q)$ be a positive function
such that $f(x_0)=0$ and $f(x) = 1$ for $x\in Q\backslash O_{x_0}$. Since
$P f \in C(Q)$ and $(Pf)(x_0) = \left<q_y, f\right> = 1$, one can
find, for an arbitrary $0 < \varepsilon < 1/2$, an open neighborhood
$V \subset O_{x_0}$ of $x_0$ such that
$f(x)<\epsilon$ and $(Pf)(x)>1-\varepsilon$ for $x\in V$.
Denote by $i_V$ the indicator of $V$. Since
$$
\left< q_{\phi(x)},  i_V f\right> = \int_{V} f(t) q_{\phi(x)} (dt)
\geq (1-\varepsilon) q_{\phi(x)}(V),
$$
we have, by using (\ref{recomp}),
\begin{eqnarray*}
\varepsilon m(V)&\geq & \int_{V} f(x) m(dx) =
\int_Q \left< q_{\phi(x)},  i_V f\right> m(dx)\\
&\geq& (1-\varepsilon)\int_{Q} q_{\phi(x)} (V) m(dx) = (1-\varepsilon) m(V).
\end{eqnarray*}
Since $m$ is positive on open sets, we have that $\ \mbox{supp}\ q_y = \phi^{-1} (y)$.
The fact that the orbital morphism $\phi$ is consistent
follows from Lemma 2.5.
\hfill$\Box$

\par {\bf 2.8. Remark.} The first condition of Theorem 2.1 is not
necessary for coassociativity of $\tilde\Delta$. Indeed, the decomposition
of ${\Bbb Z}_6 = \{0\}\cup\{1, 2\}\cup\{3\}\cup\{4, 5\}$ on $\phi$-orbits
is a DJS-hypergroup, but $\tilde\Delta(\delta_1-\delta_2) \neq 0 $, although
$\delta_1-\delta_2 \in \ker P$.

\subsection*{3. Quantum double cosets and quantum Delsart hypergroups}
\setcounter{subsection}{3}
\setcounter{equation}{0}

In this section we use Theorem 2.1 to introduce double cosets of
quantum groups \cite{ChaVai:double}, \cite{Vai:double}
and an analogue of Delsart construction for a quantum hypergroup.

\par {\bf 3.1. Quantum double cosets.} Let ${\cal A}_i = (A_i, \delta_i, \epsilon_i, \kappa_i)\ $,
$i=1,2$, be two compact matrix pseudogroups and let $\pi: A_1\to A_2$
be a Hopf $C^*$-algebra epimorphism, i. e. $\pi$ is a $C^*$-algebra
epimorphism satisfying $(\pi\otimes\pi)\circ\Delta_1 = \Delta_2\circ\pi$,
$\epsilon_2\circ\pi = \epsilon_1$ and also $\pi({A_1}_0) \subset {A_2}_0 $
with $\pi\circ\kappa_1 = \kappa_2\circ\pi$, where ${A_i}_0$ is the
$\ast$-subalgebra of $A_i$ generated by matrix elements of the fundamental
corepresentation. It was established in \cite{ChaVai} that $\pi\circ\star =
\star\circ\pi$, where $\star$ is defined by (\ref{def:star}), on each $A_i$,
$i=1, 2$. Let $\mu_i$ be a Haar measure of ${\cal A}_i$, $i=1, 2$.
The algebra $A_1$ possesses a structure of a left (right) comodule
with respect to the coactions
$\Delta_l (a) = (\pi\otimes \ \mbox{id})\Delta (a)$
(resp., $\Delta_r (a) = (\ \mbox{id}\otimes \pi)\Delta (a)$).
Define
\begin{equation}\label{def:inv}
	\begin{array}{c}
		A_1/A_2=\{a\in A_1:(id\otimes \pi)\circ
		\delta(a)=a\otimes 1\},\\[2mm]
		A_2\backslash A_1=\{a\in A_1:(\pi\otimes id)\circ
		\delta(a)=1\otimes a\},\\[2mm]
		A_2\backslash A_1/A_2=A_2\backslash A_1 \bigcap
		A_1/A_2.
	\end{array}
\end{equation}
It is immediate that $A_1/A_2$, $A_2\backslash A_1$, $A_2\backslash
A_1/A_2$ are involutive algebras with the unit $1$. Denote by
$A^{\rm inv}$ the $C^*$-algebra completion of $A_2\backslash A_1/A_2$.

\par  Define $P = \pi^l\circ\pi^r$, where
$\pi^l = (\mu_2\circ\pi\otimes \ \mbox{id})\circ\Delta_1$,
$\pi^r = (\ \mbox{id}\otimes \mu_2\circ\pi)\circ\Delta_1$
are commuting projections on $A_1$. Then {\sl $P: A_1\to A^{\rm inv}$
is a conditional expectation on $A_1$ satisfying hypotheses of
Theorem 2.1.
Hence
${\cal B} = (A^{\rm inv}, \tilde\Delta, \varepsilon_1, \star ,\sigma_t, \mu_1)$
is a compact quantum hypergroup, where $\sigma_t$ and $\star$
are defined by (\ref{def:star}), (\ref{Wor:sigma}).
Also, for all $b\in A^{\rm inv}$,
the following formula for the comultiplication $\tilde\Delta$ holds:}
\begin{equation}\label{GP:delta}
\tilde\Delta (b) = (\ \mbox{id}\ \otimes \mu_2\circ\pi\otimes\ \mbox{id}\ )
(\Delta_1\otimes\ \mbox{id})\Delta_1 (b).
\end{equation}

\par Indeed, the projections $\pi^l, \pi^r$
satisfy the equalities
\begin{equation}\label{pi:lr}
(\ \mbox{id}\ \otimes \pi^r)\circ\Delta = \Delta\circ\pi^r,\ \
(\pi^l\otimes \ \mbox{id}\ )\circ\Delta = \Delta\circ\pi^l.
\end{equation}
By using these equalities and straightforward calculation,
one can check that $\ker P$ is a coideal. Formula
(\ref{GP:delta}) easily follows from (\ref{pi:lr}):
\begin{eqnarray*}
(P\otimes P)\Delta_1 (b) &=& (\pi^r\otimes\pi^l)\circ\Delta_1(\pi^l\circ\pi^r(b))
= (\pi^r\otimes\pi^l)\circ\Delta_1(b)\\
&=& (\pi^r\otimes\ \mbox{id}\otimes\ \mbox{id})\circ(\ \mbox{id}\otimes\mu_2
\circ\pi\otimes\ \mbox{id})\circ(\ \mbox{id}\otimes\Delta_1)\circ\Delta_1(b)\\
&=& (\pi^r\otimes\ \mbox{id}\otimes\ \mbox{id})\circ(\ \mbox{id}\otimes\mu_2
\circ\pi\otimes\ \mbox{id})\circ(\Delta_1\otimes\ \mbox{id})\circ\Delta_1(b)\\
&=& (\pi^r\otimes\ \mbox{id})\circ\Delta_1 (b)
= (\ \mbox{id}\ \otimes \mu_2\circ\pi\otimes\ \mbox{id}\ )\circ
(\Delta_1\otimes\ \mbox{id})\circ\Delta_1 (b).
\end{eqnarray*}
It is obvious that $P$ is $\mu_1$-invariant:
\begin{eqnarray*}
\mu\circ P(a)  &=& (\mu_1\otimes\mu_2\circ\pi)\circ\Delta_1\circ(\mu_2\circ\pi\otimes\ \mbox{id})\circ\Delta_1(a)\\
&=& (\mu_2\circ\pi\otimes\mu_1\otimes\mu_2\circ\pi)\circ(\ \mbox{id}\otimes\Delta_1)\circ\Delta_1(a)\\
&=& \mu_1(a) (\mu_2\circ\pi\otimes\mu_2\circ\pi)(1\otimes 1) = \mu_1(a).
\end{eqnarray*}
It is also obvious that $P$ commutes with $\kappa_1$ and
that condition (A) of Remark 2.2 holds.
In virtue of Remark 2.2 we only need to examine
Condition 4 of Theorem 2.1.
To prove Condition 4, we have by using (\ref{GP:delta}) that
\begin{eqnarray*}
(\varepsilon_1\otimes\ \mbox{id}\ )\circ\tilde\Delta(a) &=&
(\varepsilon_1\otimes\ \mbox{id}\ )\circ(\ \mbox{id}\otimes\mu_2\circ\pi\otimes\ \mbox{id})
\circ(\Delta_1\otimes\ \mbox{id}\ )\circ\Delta_1(a)\\
&=& (\mu_2\circ\pi\otimes\ \mbox{id})\circ(\varepsilon_1\otimes\ \mbox{id}\otimes\ \mbox{id})
\circ(\Delta_1\otimes\ \mbox{id})\circ\Delta_1 (b)\\
&=& (\mu_2\circ\pi\otimes\ \mbox{id})\Delta_1(b)
= \pi^l (b) = b
\end{eqnarray*}
for all $b\in A_2\backslash A_1/A_2$. Similarly,
$(\ \mbox{id}\otimes\varepsilon_1)\circ\tilde\Delta(b) = b$.
\hfill $\Box$

\par {\bf 3.2. Remark.} If quantum groups ${\cal A}_1, \ {\cal A}_2$ are finite
dimensional, then {\sl condition 4 of Theorem 2.1 follows from
Proposition 2.4.} Indeed, it is obvious that $\pi^l$ (resp. $\pi^r$)
is a conditional expectation from $A_1$ to $A_2\backslash A_1$
(resp. to $A_1/A_2$). Conditions (\ref{p_epsilon1}), (\ref{p_epsilon2})
of Proposition 2.4 then follow from the following relations
\begin{eqnarray*}
(P\otimes P)((a\otimes 1)(\tilde\Delta(b))) &=& (P\otimes P)((Pa\otimes 1)(\tilde\Delta(b)))\\
(P\otimes P)((1\otimes a)(\tilde\Delta(b))) &=& (P\otimes P)((1\otimes Pa)(\tilde\Delta(b))),
\end{eqnarray*}
which hold for all $b\in A_2\backslash A_1/A_2,\ a\in A_1$.
The last relations follows by direct computations from the fact
that $\pi^l$ and $\pi^r$ are conditional expectations.

\par {\bf 3.3. Delsart hypergroups.} Let ${\cal A} = (A, \delta, \epsilon, \kappa)\ $
be a compact matrix pseudogroup and let $\Gamma$
be a compact group of Hopf $C^*$-algebra automorphisms of ${\cal A}$,
i. e. each $\gamma\in\Gamma$ is a $C^*$-algebra
automorphism satisfying $(\gamma\otimes\gamma)\circ\Delta = \Delta\circ\gamma$,
$\epsilon\circ\gamma = \epsilon$ and also $\gamma(A_0) \subset A_0 $
with $\gamma\circ\kappa = \kappa\circ\gamma$, where $A_0$ is the
$\ast$-subalgebra of $A$ generated by matrix elements of the fundamental
corepresentation. Let $\nu$ be a Haar measure of $\Gamma$ such that
$\nu(\Gamma) = 1$. In what follows, we denote integration with respect to $\nu$ by
$d\gamma$. Denote by $B=\{b\in A| \gamma(b)=b \}$ the fixed point algebra
for the $\Gamma$-action. Define a mapping $P: A\to B$ by
$P a = \int_\Gamma \gamma (a)$. Then {\sl $P$ is a conditional expectation
on $A$ satisfying hypothesis of Theorem 2.1. Hence,
${\cal B} = (B, \tilde\Delta, \varepsilon, \star ,\sigma_t,\mu)$
is a compact quantum hypergroup (Delsart hypergroup),
where $\sigma_t$ and $\star$ are defined by
(\ref{def:star}), (\ref{Wor:sigma}). Also, for all $b\in B$,
the following formula for the comultiplication $\tilde\Delta$ holds:}
\begin{equation}\label{D:delta}
\tilde\Delta (b) = (\ \mbox{id}\ \otimes P)(\Delta(b).
\end{equation}

\par Indeed, it is obvious that $P$ is a conditional expectation.
Since, for any $\gamma\in \Gamma$ and $a\in A$,
$$
(\mu\circ\gamma)(a) 1 = (\ \mbox{id}\otimes\mu)\Delta(\gamma(a)) =
\gamma((\ \mbox{id}\otimes\mu\circ\gamma)\Delta(a)),
$$
we have that $\mu\circ\gamma$ is a Haar measure. Since a normalized Haar
measure is unique, we have that $\mu\circ\gamma = \mu$, whence
it follows that $P$ is $\mu$-invariant. Since
$\varepsilon\circ\gamma = \varepsilon$ we get that
$\varepsilon\circ P = \varepsilon$ and condition 4 of
Theorem 2.1 follows from Remark 2.3.
The facts that $\ker P$ is a coideal
and that $P$ commutes with $\kappa$ are obvious.
Thus the statement follows from Remark 2.2.
\hfill $\Box$

\par {\bf 3.4.} In order to construct nontrivial examples
of finite dimensional quantum Delsart hypergroups, we need
to know about automorphisms of nontrivial Kac algebras.
A number of examples of nontrivial finite dimensional
Kac algebras are constructed in \cite{Vai:twist}, \cite{Nik:twist}
as twistings of the Kac algebras of finite groups.
It is natural to expect that automorphisms of twisted Kac
algebras are related with automorphisms of the corresponding
finite groups. In what follows, we state some results in this direction.

 Let $({\Bbb C}(G), \Delta_\Omega, \varepsilon, \kappa_\Omega)$
be a Kac algebra obtained by twisting from the cocommutative Kac algebra
$({\Bbb C}(G), \Delta, \varepsilon, \kappa)$ of a finite group $G$ with
respect to an abelian subgroup $H$ of $G$. Let  $\alpha$ be an automorphism
of $G$. Define $\alpha(\lambda(g)) = \lambda(\alpha(g))$. Then
$\alpha : {\Bbb C}(G)\to {\Bbb C}(G) $ is an automorphism of the Kac algebra
$({\Bbb C}(G), \Delta, \varepsilon, \kappa)$. Denote by
$\ \mbox{Ad} \ \Omega$ the  automorphism of ${\Bbb C}(G)\otimes{\Bbb C}(G)$
given by $\ \mbox{Ad}\  \Omega (x) = \Omega x \Omega^*$,
where $x\in {\Bbb C}(G)\otimes{\Bbb C}(G)$.

\par The next two propositions
give sufficient conditions for $\alpha$ (or a certain automorphism
constructed from $\alpha$) to be an automorphism of the
twisted Kac algebra
$({\Bbb C}(G), \Delta_\Omega, \varepsilon, \kappa_\Omega)$.

\par {\bf 3.5. Proposition.} {\sl Let $\Omega$ be a
(pseudo)-coinvolutive 2-(pseudo)-cocycle of ${\Bbb C}(G)$.
If $\alpha\otimes\alpha$
commutes with $\ \mbox{Ad} \ \Omega$, then $\alpha$ is an automorphism
of the twisted Kac algebra
$({\Bbb C}(G), \Delta_\Omega, \varepsilon, \kappa_\Omega)$.
In particular, it is sufficient  that
$\alpha\upharpoonright H = \ \mbox{id}$. }

\par {\bf Proof.} The first statement is obvious. The last
statement follows from the fact that
$\Omega \in {\Bbb} C (H)\otimes {\Bbb C} (H)$.
\hfill$\Box$

\par {\bf 3.6. Proposition.} {\sl Let $\Omega$ be a
(pseudo)-coinvolutive 2-(pseudo)-cocycle of ${\Bbb C}(G)$ and
$u = m(\ \mbox{id}\otimes\kappa)\Omega$ be the corresponding
unitary in ${\Bbb C}(G)$. Denote by $\gamma$ the automorphism
of this $C^*$-algebra defined by $\gamma(x) = u\alpha(x)u^*$
for $x\in {\Bbb C}(G)$. Suppose that the following conditions hold:
\begin{itemize}
   \item[1.] $\alpha (u) = u^*$,
   \item[2.] the element $(\alpha^{-1}\otimes\alpha^{-1})(\Omega^u)\Omega^* $
             belongs to the commutant of $\Delta ({\Bbb C}(G))$,
   \item[3.]  $\varepsilon(u) = 1$.
\end{itemize}
Then $\gamma$ is an automorphism of the coalgebra
$({\Bbb C}(G), \Delta_\omega, \varepsilon)$, i. e.
$(\gamma\otimes\gamma)\circ\Delta_\Omega = \Delta_\Omega\circ\gamma$,
$\varepsilon\circ\gamma = \varepsilon$.
If $\kappa (u) = u$, then $\gamma$ is an automorphism of the
Kac algebra $({\Bbb C}(G), \Delta_\omega, \varepsilon, \kappa_\Omega)$.
The last condition is always true when $\Omega$ is a counital 2-cocycle,
but if $\Omega$ is only a 2-pseudococycle, then one should verify
that $\gamma$ commutes with $\kappa_\Omega$.
}

\par {\bf Proof.} It follows from the first condition of the proposition
that $\gamma = \ \mbox{Ad}\ u\circ\alpha = \alpha\circ\ \mbox{Ad}\ u^*$.
Following condition 2 one can find an element $Z$ in commutant of
$\Delta ({\Bbb C}(G))$ such that
$(\alpha^{-1}\otimes\alpha^{-1})(\Omega^u) = \Omega  Z$.
Then, for any $x\in {\Bbb C}(G)$, we have
\begin{eqnarray*}
\Delta_\Omega (\gamma(x)) &=& \Omega\Delta (u\alpha(x)u^*)\Omega^*
= \Omega\Delta (u)((\alpha\otimes\alpha)\Delta(x))\Delta (u^*)\Omega^* \\
&=& (\alpha\otimes\alpha)
\left((\alpha^{-1}\otimes\alpha^{-1})(\Omega\Delta(u)) \Delta(x)
(\alpha^{-1}\otimes\alpha^{-1})(\Delta(u^*)\Omega^*)\right)\\
&=& (\gamma\otimes\gamma)
\left((u\otimes u)(\alpha^{-1}\otimes\alpha^{-1})(\Omega\Delta(u)) \Delta(x)
(\alpha^{-1}\otimes\alpha^{-1})(\Delta(u^*)\Omega^*)(u^*\otimes u^*)\right)\\
&=& (\gamma\otimes\gamma)
\left((\alpha^{-1}\otimes\alpha^{-1})((u^*\otimes u^*)\Omega\Delta(u)) \Delta(x)
(\alpha^{-1}\otimes\alpha^{-1})(\Delta(u^*)\Omega^*(u\otimes u))\right)\\
&=& (\gamma\otimes\gamma)
\left((\alpha^{-1}\otimes\alpha^{-1})(\Omega^u) \Delta(x)
(\alpha^{-1}\otimes\alpha^{-1})(\Omega^u)^*\right)\\
&=& (\gamma\otimes\gamma)(\Omega\Delta(x)\Omega^*) \\
&=& (\gamma\otimes\gamma)(\Delta_\Omega (x)).
\end{eqnarray*}
It follows from condition 3 of the proposition  that
$\varepsilon\circ\gamma = \varepsilon$.
To prove the last statement, let us suppose that $\kappa (u) = u$.
Then, for any $x\in {\Bbb C}(G)$, we have $\gamma(\kappa_\Omega (x))
= u\alpha (u\kappa (x) u^*)u^* = \alpha(\kappa(x)) =
\kappa(\alpha(x))$. On the other hand,
$\kappa_\Omega (\gamma(x)) = u\kappa(u\alpha(x)u^*)u^*
= u\kappa(u^*)\kappa(\alpha(x))\kappa(u)u^* = \kappa(\alpha(x)).$
\hfill$\Box$

\par {\bf 3.7. Remark.} By using the action of an automorphism $\alpha$
on the group $\hat H$ which is dual to the Abelian subgroup $H$,
$<\alpha (\hat h), h> = < \hat h, \alpha^{-1}(h)>$,
one can rewrite condition 1 of Proposition 3.7 in
terms of the 2-(pseudo)-cocycle $\omega$ on $\hat H\times\hat H$,
$$
\omega(\alpha(\hat x), \alpha(\hat x^{-1})) = \overline{\omega(\hat x, \hat x^{-1})} \quad (\hat x\in \hat H).
$$
Condition 2 of Proposition 3.7 follows from the relation
$(\alpha^{-1}\otimes\alpha^{-1})(\Omega^u) = \Omega$ which
can also be rewritten in terms of the 2-(pseudo)-cocycle $\omega$,
$$
\omega(\hat x, \hat x^{-1})\omega(\hat y, \hat y^{-1})\omega(\alpha(\hat x), \alpha(\hat y))
= \omega (\hat x, \hat y)\omega(\hat x\hat y, (\hat x\hat y)^{-1})\quad
(\hat x, \hat y\in \hat H).
$$
One can also rewrite other equalities that contain $\Omega$ in such manner,
for example, the relation $(\alpha^{-1}\otimes\alpha^{-1})\Omega = \Omega$ is sufficient
for the statement of Proposition 3.6 to hold,
\begin{equation}\label{omega}
\omega(\alpha(\hat x), \alpha(\hat y)) = \overline{\omega(\hat x, \hat y)}.
\end{equation}

\subsection*{4. Examples}
\setcounter{subsection}{4}
\setcounter{equation}{0}

\par{\bf 4.1. Quantum hypergroups associated with a twisting of the quasiquaternionic group.}
Let $G = Q_n$ ($n = 2, 3, \dots $) be the quasiquaternionic group generated by two
elements, $a$ of order $2n$ and $b$ of order 4 such that
$b^2 = a^n$ and $bab^{-1} = a^{-1}$. The group
$Q_n = \{a^k, ba^k, k = 0, 1, \dots , 2n-1 \}$
and the group algebra $A = {\Bbb C}(G)$ is isomorphic to
$${\Bbb C}\oplus{\Bbb C}\oplus{\Bbb C}\oplus{\Bbb C}
\oplus \underbrace{M_2({\Bbb C})\oplus\dots\oplus
M_2({\Bbb C})}_{n-1}.$$
Let $e_1$, $e_2$, $e_3$, $e_4, $$e^{j}_{11}$, $e^{j}_{12}$,
$e^{j}_{21}$,
$e^{j}_{22}\ (j=1,...,n-1)$, be the matrix units of this algebra. We
can now write the left regular representation $\lambda$ of $G$ \cite{HR},
$$\lambda (a^k)=e_1 +e_2 +(-1)^k(e_3+ e_4) +\sum_j
(\varepsilon^{jk}_n
e^{j}_{11}+\varepsilon^{-jk}_n e^{j}_{22}),$$
$$\lambda (a^k b)=e_1 -e_2 + (-1)^k(e_3- e_4) +\sum_j
(\varepsilon^{j(k-n)}_n e^{j}_{12}+\varepsilon^{-jk}_n e^{j}_{21}),$$
for even $n$, and
$$\lambda (a^k b)=e_1 -e_2 + i(-1)^k(e_3- e_4) +\sum_j
(\varepsilon^{j(k-n)}_n e^{j}_{12}+\varepsilon^{-jk}_n e^{j}_{21}),$$
for odd $n$, where $\varepsilon_n=e^{i\pi/n}$.

\par Consider the subgroup $H={\Bbb Z}_4=\{e,b,b^2,b^3\}$. Since the dual
group $\hat H$ is isomorphic to $H$, following \cite{Vai:twist},
one can compute the orthogonal projections $P_{\hat h}\ \ (\hat h\in \hat H)$,
\begin{eqnarray*}
P_{\hat e} &=& e_1+e_3+q_1,\ P_{\hat b} =p_1, \\
P_{{\hat b}^2}&=&e_2+e_4+q_2,\ P_{{\hat b}^3}=p_2,\ for\ even\ n,\\
P_{\hat e}&=&e_1+q_1,\ P_{\hat b} =e_4+p_1,\\
P_{{\hat b}^2}&=&e_2+q_2,\ P_{{\hat b}^3}=e_3+p_2,\ for\ odd\ n,
\end{eqnarray*}
where $p_1$, $p_2$, $q_1$, $q_2$ are the orthogonal projections
defined by :
\begin{eqnarray*}
p_1&=&1/2 {\sum}'_j(e^{j}_{11}-ie^{j}_{12}+ie^{j}_{21}+
e^{j}_{22}),\\
p_2&=&1/2 {\sum}'_j(e^{j}_{11}+ie^{j}_{12}-ie^{j}_{21}+
e^{j}_{22}),\\
q_1&=&1/2 {\sum}''_j(e^{j}_{11}+e^{j}_{12}+e^{j}_{21}+ e^{j}_{22}),\\
q_2&=&1/2 {\sum}''_j(e^{j}_{11}-e^{j}_{12}-e^{j}_{21}+ e^{j}_{22}),
\end{eqnarray*}
where ${\sum}'$ (resp., ${\sum}''$) means that the corresponding
index in the summation takes only odd (resp., even) values.
One can see that $e_1$ is the projection given by the co-unit of
$A$.

\par The Abelian subalgebra generated by
$\lambda(e)$, $\lambda(b)$,$\lambda(b^2)$, $\lambda(b^3)$
is also generated by the mutually orthogonal orthogonal projections
$P_{\hat e}, \ P_{\hat b},\  P_{{\hat b}^2}, \ P_{{\hat b}^3}$.
The projections $P_{\hat e}+P_{{\hat b}^2}$ and
$P_{\hat b}+ P_{{\hat b}^3}$ are central.

\par The unitary $\Omega\in A\otimes A$
is obtained by the lifting construction
from the pseudo-cocycle $\omega$ on
$\hat H\times\hat H$ such that,
for all $u$, $h$ in $\hat H$,
$$\omega({\hat e},u)=\omega(u,u)=1, \overline{\omega(h,u)}=
\omega (u,h)$$
and $\omega({\hat b}, {\hat b}^2)=\omega({\hat b}^2, {\hat b}^3)=
\omega({\hat b}^3,{\hat b})=i$.
In \cite{Vai:twist} it was established that
\begin{eqnarray*}
\Omega &=& P_{\hat e}\otimes I+P_{\hat b}\otimes(P_{\hat e}+P_{\hat
b} + i(P_{{\hat b}^2}-P_{{\hat b}^3}))+P_{{\hat b}^2}\otimes(P_{\hat e}+
P_{{\hat b}^2}\\
&+& i(P_{{\hat b}^3}-P_{\hat b}))+P_{{\hat
b}^3}\otimes
(P_{\hat e}+P_{{\hat b}^3}++i(P_{\hat b}-P_{{\hat b}^2}))
\end{eqnarray*}
is a pseudo-coinvolutive 2-pseudo-cocycle with respect to the unitary
$$
u = P_{\hat e} + P_{{\hat b}^2} + i(P_{{\hat b}^3} - P_{\hat b})
$$
such that ${\cal A}_\Omega = (A, \Delta_{\Omega}, \varepsilon, \kappa_{\Omega}, \mu)$
is a nonsymmetric Kac algebra
(i. e. $\Sigma\circ\Delta_{\Omega}\neq\Delta_{\Omega}$).

\par Define an authomorphism $\alpha$ of order 2 on $G=Q_n$:
$\alpha(a) = a$, $\alpha(b) = b^3$.
Let us verify that $\alpha$ satisfyes the hypothesis of
Proposition 3.7. Indeed, $\alpha (P_{\hat b}) =
P_{{\hat b}^3}$ and, hence, $\alpha (u) = u^*$.
Let $Z = (\alpha\otimes\alpha)(\Omega^u) \Omega^*$.
Then it is not hard to compute that
$$
Z = (P_{\hat e}+ P_{{\hat b}^2})\otimes I +
(P_{\hat b} + P_{{\hat b}^3})\otimes
((P_{\hat e}+ P_{{\hat b}^2}) - (P_{\hat b} + P_{{\hat b}^3})) \in
Z(A)\otimes Z(A)\subset \Delta(A)',
$$
where $Z(A)$ is the center
of $A$. In order to show that $\gamma(x) = u\alpha(x)u^*$,
$x\in {\Bbb C}(G)$, is an authomorphism of the Kac algebra
${\cal A}_\Omega $, we need to verify that
$\gamma\circ\kappa_\Omega = \kappa_\Omega\circ\gamma$.
Indeed, it is obvious that $\kappa (u) = u^*$. Then
$\kappa_\Omega (\gamma(x)) = u^2 \alpha(\kappa(x)) (u^*)^2$ and
$\gamma(\kappa_\Omega(x)) = \alpha(\kappa(x))$ for all $x\in A$. But
$$
u^2 = P_{\hat e} + - P_{\hat b} + P_{{\hat b}^2} - (P_{{\hat b}^3} = \lambda(b^2)\in Z(A).
$$
Hence, $\gamma$ is an authomorphism of the Kac algebra ${\cal A}_\Omega$ of
order 2.

\par Let us consider the group $\Gamma = \{\mbox{id}, \gamma\}$ of authomorphisms
of the Kac algebra ${\cal A}_\Omega$ and define a conditional expectation
$P$ on ${\cal A}_\Omega$ via the formula $P (x) = \frac 1 2 (x + \gamma (x))$.
Thus we obtain a Delsart hypergroup $(B = P(A), \tilde\Delta_\Omega, \varepsilon,
\kappa_\Omega, \mu)$.
In order to determine the structure of $B$, note that
$\alpha(e_3) = e_4$ and $\alpha(e^j_{12}) = (-1)^j e^j_{12}\ (j = 1, \dots , n-1)$
for odd $n$ and $\alpha (e_i) = e_i\ (i=1,\dots , 4)$,
$\alpha(e^j_{12}) = (-1)^j e^j_{12}$, $\alpha(e^j_{21}) = (-1)^j e^j_{21}
\ (j = 1, \dots , n-1)$ for even $n$.
Since, for odd $n$,
$$
u = e_1 + e_2 + i (e_3 - e_4) + {\sum_j}' (e^j_{21} - e^j_{12}) +
{\sum_j}'' (e^j_{11} - e^j_{22})
$$
and, for even $n$,
$$
u = e_1 + e_2 + e_3 + e_4 + {\sum_j}' (e^j_{21} - e^j_{12}) +
{\sum_j}'' (e^j_{11} - e^j_{22}),
$$
we obtain an explicit formula for the action of $\gamma$. Indeed, for odd $n$,
$ \gamma(e_3) = e_4 $, for even $j$, we have $\gamma(e^j_{kl}) = e^j_{kl}$,
$k, l = 1, 2$,
and for odd $j$, we have $\gamma(e^j_{11}) = e^j_{22}$,
$\gamma(e^j_{12}) = e^j_{21}$ . For even $n$, we have
$\gamma(e_k) = e_k$, $k = 1, \dots , 4$,
$\gamma(e^j_{kl}) = e^j_{kl}$, $k, l = 1, 2$, for even $j$, and
$\gamma(e^j_{11}) = e^j_{22}$,
$\gamma(e^j_{12}) = e^j_{21}$ for odd $j$.
Hence for odd $n$, we have
$$
B_n = \underbrace{{\Bbb C}\oplus\dots\oplus{\Bbb C}}_{n+2}\oplus
\underbrace{M_2({\Bbb C})\oplus\dots\oplus
M_2({\Bbb C})}_{\frac {n-1} 2}
$$
and, for even $n$,
$$
B_n = \underbrace{{\Bbb C}\oplus\dots\oplus{\Bbb C}}_{n+4}\oplus
\underbrace{M_2({\Bbb C})\oplus\dots\oplus
M_2({\Bbb C})}_{\frac {n-2} 2}.
$$
In any case, $\dim B = 3n$.

\par Let us show that the comultiplication $\tilde\Delta_\Omega$ is not
symmetric. For this, it is enough to prove that
\begin{equation}\label{nococo}
(P\otimes P)(\Delta_\Omega)(\lambda(ab)) - \Sigma\Delta_\Omega (\lambda(ab)))\neq 0,
\end{equation}
where $\Sigma$ is the flip in $A\otimes A$.
It is useful to note that $\Omega = \Omega_1 + i\Omega_2$,
where $\Sigma\Omega_1 = \Omega_1^* = \Omega_1$ and
$\Sigma\Omega_2 = -\Omega_2^* = -\Omega_2$. Then relation
(\ref{nococo}) is equivalent to
$$
(P\otimes P)(\Omega_1(\lambda(ab)\otimes\lambda(ab))\Omega_2 - \Omega_2(\lambda(ab)\otimes\lambda(ab))\Omega_1) \neq 0.
$$
By direct calculations one can obtain that
\begin{eqnarray*}
(P\otimes P)(\Omega_1(\lambda(ab)\otimes\lambda(ab))\Omega_2 - \Omega_2(\lambda(ab)\otimes\lambda(ab))\Omega_1) \\
= {\sum_k}' \cos \frac {\pi k} n (e^k_{11} + e^k_{22})\otimes {\sum_l}{''} \sin \frac {\pi l} n (e^l_{12} - e^l_{21})\neq 0
\end{eqnarray*}
for even $n$ and
\begin{eqnarray*}
(P\otimes P)(\Omega_1(\lambda(ab)\otimes\lambda(ab))\Omega_2 - \Omega_2(\lambda(ab)\otimes\lambda(ab))\Omega_1) \\
= {\sum_k}' \cos \frac {\pi k} n (e^k_{11} + e^k_{22})\otimes {\sum_l}{''} \sin \frac {\pi l} n (-e^l_{11} + e^l_{22})\neq 0
\end{eqnarray*}
for odd $n$, where $\Sigma'$ (resp. ${\Sigma}{''}$) means that the corresponding index in the
summation takes only odd (resp. even) values.

\par {\bf 4.2. Remark.} The Kac algebra ${\cal A}_\Omega$ for $n=2$
is nothing else but the historical Kac-Paljutkin example of a non-trivial
Kac algebra \cite{KP}. In this
case, the algebra $B_2$ is commutative but the comultiplication $\tilde\Delta_\Omega$
is not symmetric. So $B_2$ is the usual noncommutative hypergroup of order $6$.
For $n\geq 3$, the quantum hypergroups ${\cal B}_n$ are nontrivial.
We obtain exact formulas for the comultiplication in
$$
B_3 = \underbrace{{\Bbb C}\oplus\dots\oplus{\Bbb C}}_5\oplus M_2({\Bbb C}).
$$
Since $Q_3$ is isomorphic to ${\Bbb Z}_3\rtimes {\Bbb Z}_4$, one can use
the left regular representation of ${\Bbb Z}_3\rtimes {\Bbb Z}_4$
instead of $Q_3$ (this means that the basis of matrix units of
the algebra $A$
is changed up to corresponding authomorphism which does not permit two-
dimensional minimal ideals).
The formulas for comultiplication in the new basis become simplier.
So, in the basis of $B_3$: $f_1 = P (e_1) = e_1, \ f_2 = P (e_2) = e_2, \
2 f_3 = P (e_3 ) = 1/2 (e_3 + e_4 ) = P (e_4 ), \
2 f_4 = P ( e^2_{11} ) = 1/2 ( e^2_{11} + e^2_{22} ) = P( e^2_{22} ), \
2 f_5 = P ( e^2_{12} ) = 1/2 ( e^2_{12} + e^2_{21} ) = P( e^2_{21} ), \
f_{ij} = P (e^1_{ij}) = e^1_{ij}$, $i, j = 1, 2$, we have:
\begin{eqnarray*}
\tilde\Delta_\Omega (f_1) &=& f_1 \otimes f_1 + f_2\otimes f_2 + 1/2 f_3 \otimes f_3
+ 1/4 f_4 \otimes f_4 + 1/4 f_5 \otimes f_5\\
&+& 1/2\left( f_{11}\otimes f_{22} + f_{12}\otimes f_{21}
+ f_{21}\otimes f_{12} + f_{22}\otimes f_{11} \right),
\\
\tilde\Delta_\Omega (f_2) &=&  f_1 \otimes f_2 + 1/2 f_3 \otimes f_3 + f_2\otimes f_1
+ 1/4 f_4 \otimes f_4 + 1/4 f_5 \otimes f_5\\
&+& 1/2\left( f_{11} \otimes f_{22} - f_{12}\otimes f_{21}
- f_{21}\otimes f_{12} + f_{22}\otimes f_{11} \right),
\\
\tilde\Delta_\Omega (f_3) &=& (f_1 + f_2) \otimes f_3 + f_3 \otimes (f_1 + f_2)
+ 1/2\{ (f_{11} + f_{22})\otimes f_4\\
&+& f_4 \otimes (f_{11} + f_{22}) + (f_{11} - f_{22})\otimes f_5 - f_5 \otimes (f_{11} - f_{22}),
\\
\tilde\Delta_\Omega (f_4) &=& (f_1 + f_2) \otimes f_4 + f_4 \otimes (f_1 + f_2)
+ (f_{11} + f_{22}) \otimes f_3 + f_3 \otimes (f_{11} + f_{22})\\
&+& 1/2\{ (f_{11} + f_{22})\otimes f_4 + f_4 \otimes (f_{11} + f_{22})
- (f_{11} - f_{22})\otimes f_5\\
&+& f_5 \otimes (f_{11} - f_{22}) \},
\\
\tilde\Delta_\Omega (f_5) &=& (f_1 + f_2) \otimes f_5 + f_5 \otimes (f_1 + f_2)
- (f_{11} - f_{22}) \otimes f_3 + f_3 \otimes (f_{11} - f_{22})\\
&+& 1/2\{ (f_{11} - f_{22})\otimes f_4 - f_4 \otimes (f_{11} - f_{22})
- (f_{11} + f_{22})\otimes f_5\\ &-& f_5 \otimes (f_{11} + f_{22})\},
\\
\tilde\Delta_\Omega (f_{11}) &=& f_1 \otimes f_{11} + f_{11} \otimes f_1
+ f_2 \otimes f_{11} + f_{11}\otimes f_2
+ 1/2 f_3 \otimes \left( f_4 + f_5 \right)\\
&+& 1/2 \left( f_4 - f_5 \right) \otimes f_3
+ f_{22}\otimes f_{22} + 1/4 (f_4 + f_5) \otimes (f_4 - f_5),
\\
\tilde\Delta_\Omega (f_{12}) &=& f_1 \otimes f_{12} + f_{12} \otimes f_1
- f_2\otimes f_{12} -  f_{12}\otimes f_2 + f_{21}\otimes f_{21},
\\
\tilde\Delta_\Omega (f_{21}) &=& f_1 \otimes f_{21} + f_{21} \otimes f_1
- f_2\otimes f_{21} -  f_{21}\otimes f_2
+ f_{12}\otimes f_{12},
\\
\tilde\Delta_\Omega (f_{22}) &=& f_1 \otimes f_{22} + f_{22} \otimes f_1
+ f_2 \otimes f_{22} + f_{22}\otimes f_2
+ 1/2 f_3 \otimes \left( f_4 - f_5 \right)\\
&+& 1/2 \left( f_4 + f_5 \right) \otimes f_3
+ f_{11}\otimes f_{11} +
1/4 (f_4 - f_5) \otimes (f_4 + f_5).
\end{eqnarray*}

\par {\bf 4.3. Quantum hypergroups associated with a twisting of the dihedral group.}
Let $G=D_{2n}={\Bbb Z}_{2n}\rtimes_\alpha{\Bbb Z}_2
\ (2\leq n\in {\Bbb N})$ be the dihedral group  with the following action
of ${\Bbb Z}_2=\{1,b\}$ on ${\Bbb Z}_{2n}=\{a^k\ (k=0,1,...,2n-1)\}$
$$
\alpha_b(a^k)=a^{2n-k}\ (k=0,1,...,2n-1).
$$
The group algebra $A = {\Bbb C}(G)$ is
isomorphic to
$${\Bbb C}\oplus{\Bbb C}\oplus{\Bbb C}\oplus{\Bbb C}
\oplus \underbrace{M_2({\Bbb C})\oplus\dots\oplus
M_2({\Bbb C})}_{n-1}.$$

\par Let $e_1$, $e_2$, $e_3$, $e_4, $$e^{j}_{11}$, $e^{j}_{12}$,
$e^{j}_{21}$,
$e^{j}_{22}\ (j=1,...,n-1)$ be the matrix units of this algebra; we
can now write the left regular representation $\lambda$ of $G$ (\cite{HR},
27.61):
\begin{eqnarray*}
\lambda (a^k) &=& e_1 +e_2 +(-1)^k(e_3+ e_4) +\sum_j
(\varepsilon^{jk}_n
e^{j}_{11}+\varepsilon^{jk(2n-1)}_n e^{j}_{22}),\\
\lambda (ba^k) &=& e_1 -e_2 -(-1)^k(e_3- e_4) +\sum_j
(\varepsilon^{jk(2n-1)}_n
e^{j}_{12}+\varepsilon^{jk}_n e^{j}_{21}),
\end{eqnarray*}
where $\varepsilon_n=e^{i\pi/n}$.

\par Consider the Abelian subgroup $H=\{e,a^n,b,ba^n\}$
which is isomorphic to ${\Bbb Z}_2\times{\Bbb Z}_2$. Since the dual
group $\hat H$ is isomorphic to $H$,  one can compute the
orthogonal projections $P_{\hat h},\ (\hat h\in \hat H)$ (see
\cite{Vai:twist}):  $$P_{\hat e} =e_1+{{1+(-1)^n}\over 2}e_4+q_1, $$
$$
P_{\hat b} =e_2+{{1+(-1)^n}\over 2}e_3+q_2,
$$
$$
P_{{\hat a}^n}={{1-(-1)^n}\over 2}e_4+p_1,
$$
$$
P_{{{\hat b}{\hat a}}^n}={{1-(-1)^n}\over 2}e_3+p_2,
$$
where $p_1$, $p_2$, $q_1$, $q_2$ are the orthogonal projections
defined by :
$$
p_1=1/2 {\sum}'_j(e^{j}_{11}+e^{j}_{12}+e^{j}_{21}+ e^{j}_{22}),
$$
$$
p_2=1/2 {\sum}'_j(e^{j}_{11}-e^{j}_{12}-e^{j}_{21}+ e^{j}_{22}),
$$
$$
q_1=1/2 {\sum}''_j(e^{j}_{11}+e^{j}_{12}+e^{j}_{21}+ e^{j}_{22}),
$$
$$
q_2=1/2 {\sum}''_j(e^{j}_{11}-e^{j}_{12}-e^{j}_{21}+ e^{j}_{22}),
$$
where ${\sum}'$ (resp., ${\sum}''$) means that the corresponding
index in the
summation takes only odd (resp., even) values.
The orthogonal projections $P_{\hat e}+P_{\hat b}$ and $P_{{\hat
a}^n} + P_{{{\hat b}{\hat a}}^n}$ are central.

\par Let us consider the 2-cocycle $\omega$ on $\hat {H}
\times\hat H$
such that, for all $u$, $h$ in $\hat H$,
$\ \omega({\hat e},u)=\omega(u,u)=1$, $\overline{\omega(h,u)}=\omega (u,h)$,
and $\omega({\hat a}^n,{\hat b})=\omega({\hat b},{{\hat b}{\hat a}}^n)=
\omega({{\hat b}{\hat a}}^n,{\hat a}^n)=i$.

Let $\Omega$ be the 2-cocycle of $A \otimes A$ obtained by
the lifting construction,
\begin{eqnarray*}
\Omega &=& P_{\hat e}\otimes I+P_{{\hat a}^n}\otimes (P_{\hat
e}+P_{{\hat a}^n}
+i(P_{\hat b}-P_{{\hat b}{\hat a}^n}))+P_{\hat b}\otimes(P_{\hat
e}+P_{\hat b}\\
&+& i(P_{{\hat b}{\hat a}^n}-P_{{\hat a}^n}))+P_{{\hat b}{\hat
a}^n}\otimes (P_{\hat e}+P_{{\hat b}{\hat a}^n}+i(P_{{\hat a}^n}-P_{\hat b})).
\end{eqnarray*}
We can also write $\Omega = \Omega_1 + i\Omega_2$,
where $\Sigma\Omega_1 = \Omega_1 = {\Omega_1}^*$,
and $\Sigma\Omega_2 = -\Omega_2 = -{\Omega_2}^* $,
where $\Sigma$ is the flip in $A\otimes A$.
It is clear that the 2-cocycle $\Omega$ is strongly
co-involutive on $(A, \Delta,\kappa)$. So
${\cal A}_\Omega = (A, \Delta_\Omega, \varepsilon, \kappa, \mu)$
is a nontrivial Kac algebra \cite{Vai:twist}.

\par Let $\gamma (a) = a^p$ be an involute automorphism of ${\Bbb Z}_{2n}$,
where $p < 2n - 1$ has no common divisors with $2n$, $n \geq 4$,
and $p^2-1=0$ modulo $2n$.
It is clear that we can extend $\gamma$ to the group $G$ by
setting $\gamma (b) = b$. Since $\gamma$ acts trivially on $H$,
we have, by Proposition 3.6, that $\gamma$ is an automorphism
of the Kac algebra ${\cal A}_\Omega$.
Let us consider the group $\Gamma = \{\mbox{id}, \gamma\}$ of authomorphisms
of the Kac algebra ${\cal A}_\Omega$ and define a conditional expectation
$P$ on ${\cal A}_\Omega$ via the formula $P (x) = \frac 1 2 (x + \gamma (x))$.
Thus we obtain a quantum Delsart hypergroup $(B = P(A), \tilde\Delta_\Omega, \varepsilon,
\kappa, \mu)$. Since $p\neq 2n-1$, we have that $B$ is noncommutative
(note that, if $\gamma^2 \neq \mbox{id}$, this is not true in general).
The dimension of $B$ equals the number of $\Gamma$-orbits, so we have $\dim B = 2n + r$, where
$r$ is the number of solutions of the equation $x(p-1) = 0$ modulo $2n$. If
$n$ is a prime number, then $\dim B = 2n + 2$.

\par For example, if $n=4$ and $\gamma(a) = a^3$, then we have
$\gamma (e_i) = e_i$, $i = 1, \dots , 4$, $\gamma(e^1_{ij}) = e^3_{ij}$,
$i, j = 1, 2$, and $\gamma(e^2_{11}) = e^2_{22}$, $\gamma(e^2_{12}) = e^2_{21}$.
Thus
$$
B = \underbrace{{\Bbb C}\oplus \dots \oplus {\Bbb C}}_6 \oplus M_2 ({\Bbb C}).
$$

\par Let us show that the comultiplication $\tilde\Delta_\Omega$ is not
symmetric. For this it is enough to prove that
$$
(P\otimes P)(\Delta_\Omega)(\lambda(a)) - \Sigma\Delta_\Omega (\lambda(a)))\neq 0,
$$
where $\Sigma$ is the flip in $A\otimes A$, or, equivalently, that
$$
(P\otimes P)(\Omega_1(\lambda(a)\otimes\lambda(a))\Omega_2 - \Omega_2(\lambda(a)\otimes\lambda(a))\Omega_1) \neq 0.
$$
Since
\begin{eqnarray*}
\Omega_1 &=& \frac 1 2 \lambda (e)\otimes \lambda (e) +
\frac 1 8 \{\lambda (e)\otimes (\lambda (e) + \lambda (a^n) + \lambda (b) + \lambda (ba^n))\\
&+& \lambda (a^n)\otimes (\lambda (e) + \lambda (a^n) - \lambda (b) - \lambda (ba^n))\\
&+& \lambda (b)\otimes (\lambda (e) - \lambda (a^n) + \lambda (b) - \lambda (ba^n))\\
&+& \lambda (ba^n)\otimes (\lambda (e) - \lambda (a^n) - \lambda (b) + \lambda (ba^n))\}
\end{eqnarray*}
and
\begin{eqnarray*}
\Omega_2 &=& \frac 1 4 \{\lambda (a^n)\otimes \lambda (b)
- \lambda (a^n)\otimes\lambda(ba^n) -  \lambda(b)\otimes \lambda(a^n)
+ \lambda(b)\otimes\lambda(ba^n)\\
&+& \lambda(ba^n)\otimes \lambda(a^n) - \lambda(ba^n)\otimes\lambda(b)\},
\end{eqnarray*}
we have
\begin{eqnarray*}
\Omega_2 \lambda(a)&\otimes&\lambda(a) \Omega_1 - \Omega_1 \lambda(a)\otimes\lambda(a) \Omega_2 \\
&=& \frac 1 4 \lambda(a^{n-1})\otimes (\lambda(a^{2n-1}) - \lambda(a^{n-1})
- \lambda(ba^{n+1}) + \lambda(ba))\\
&+& \frac 1 2 (\lambda(a^{n+1}) - \lambda(a^{2n-1}))\otimes
\lambda(b)[(\lambda(a) - \lambda(a^{n+1}))  + (\lambda(a^{n-1}) - \lambda(a^{2n-1}))]\\
&+& \frac 1 2 \lambda(b)(\lambda (a^{n+1}) +
\lambda (a^{2n-1}))\otimes(\lambda(a^{n-1}) - a^{2n-1})\\
&+& \frac 1 2 \lambda(b)(\lambda(a) - \lambda(a^{n-1}))\otimes\lambda(b)(\lambda(a^{n+1}) - \lambda(a^{2n-1}))\\
&+& \frac 1 2 \lambda(b)(\lambda(a) + \lambda(a^{n-1}))\otimes (\lambda(a^{2n-1}) - \lambda(a^{n+1}))\\
&+& \frac 1 2 \lambda(b)(\lambda(a^{2n-1}) - \lambda(a^{n+1}))\otimes
\lambda(b)(\lambda(a) - \lambda(a^{n-1})).
\end{eqnarray*}
Apply $P\otimes P$ to the both sides of this equality. Then
the sum of the first three summands is equal to zero iff $p=n+1$,
while the sum of the second three summands is equal to zero
iff $p = n-1$. Since $n\geq 4$, we obtain that
the comultiplication $\tilde\Delta_\Omega$ is not symmetric.

\par {\bf 4.4. Quantum hypergroups associated with a twisting of the
symmetric group.} The twisting of the symmetric group $S_n$, $n\geq 4$, by a
2-cocycle lifted from the Abelian subgroup $H\cong {\Bbb Z}_2\times{\Bbb Z}_2$
generated by the permutations $a = (12)$ and $b = (34)$
was constructed in \cite{Nik:twist}. The 2-cocycle $\omega$
on $\hat H\times\hat H$ is defined by $\omega(\hat a, \hat b) =
\omega (\hat b, \hat a\hat b) = \omega (\hat a\hat b, a) = i$,
$\omega(\hat e, \hat x) = \omega(\hat x, \hat e) = \omega(\hat x, \hat x) = 1$,
$\omega(\hat x, \hat y) = \overline{\omega(\hat y, \hat x)}$
for all $\hat x, \hat y \in \hat H$. Let $\Omega$ be the lifted
counital 2-cocycle, as in 3.4.
Then the twisted Kac algebra ${\cal A}_\Omega$ is non-symmetric.

\par Denote by $\gamma$ the inner automorphism of the group $S_n$
generating by the permutation $a$. Since $\gamma$ acts trivially
on the subgroup $H$,  by Proposition 3.6 we have that
$\gamma$ is an automorphism of the twisted Kac algebra ${\cal A}_\Omega$
of order 2. Let $P$ be the conditional expectation on ${\cal A}_\Omega$
associated with the subgroup $\Gamma = \{\mbox{id}, \gamma\}$.
Then $B = P({\Bbb C}(S_n))$ is a quantum hypergroup.
The algebra $B$ is noncommutative, since $\lambda (b)$ does not
commute with $P(\lambda(c))$, where $c = (2341)$.
The dimension of $B$ equals the number of $\Gamma$-orbits, so
we have $\dim B = \frac 1 2 (n^2 - n +2)(n-2)!$.

\par Let us show that the comuliplication $\tilde\Delta_\Omega$ is not symmetric.
We can write $\Omega = \Omega_1 + i\Omega_2$,
where $\Sigma\Omega_1 = \Omega_1 = {\Omega_1}^*$
and $\Sigma\Omega_2 = -\Omega_2 = -{\Omega_2}^* $,
where $\Sigma$ is the flip in $A\otimes A$.
Thus, as in the previous example, our statement follows from the inequality
\begin{equation}\label{Sn:noco}
(P\otimes P)(\Omega_1(\lambda(c)\otimes\lambda(c))\Omega_2 - \Omega_2(\lambda(c)\otimes\lambda(c))\Omega_1) \neq 0
\end{equation}
which can be obtained by strightforward calculations with $c=(2341)$.

\par {\bf 4.5. Quantum hypergroups associated with a twisting of the
alternating group.}
A nontrivial twisting of the alternating group $A_n$,
$n\geq 4$, by a 2-cocycle lifted from the abelian subgroup $H\cong {\Bbb
Z}_2\times{\Bbb Z}_2$ generated by the elements $a = (12)(34)$ and $b =
(13)(24)$ was constructed in \cite{Nik:twist}. The 2-cocycle $\omega$ on
$\hat H\times\hat H$ is the same as in 4.5. Let $\Omega$ be the lifted
counital 2-cocycle, as in 3.4.  Then the twisted Kac algebra ${\cal
A}_\Omega$ is non-symmetric iff $n\geq 5$.

\par Let $\gamma$ be the restriction  of the inner
automorphism of the group $S_n$ generated by the permutation $(12)$
to the group $A_n$.
It is easy to verify that relation (\ref{omega}) holds for
the choosen automorphism $\gamma$. Hence, $\gamma$ is an automorphism
of the twisted Kac algebra ${\cal A}_\Omega$ by virtue of Remark 3.8.
Let $P$ be the conditional expectation on ${\cal A}_\Omega$
associated with the subgroup $\Gamma = \{\mbox{id}, \gamma\}$.
Then $B = P({\Bbb C}(A_n))$ is the quantum hypergroup.
It is obvious that the algebra $B$ is noncommutative
and $\dim B = \frac 1 4 (n^2 - n + 2)(n-2)!$.

\par We can obtain, by strightforward calculations,
that ineqiality (\ref{Sn:noco}) holds with
$c = (345)\in A_n$, $n\geq 5$. Thus
the comuliplication $\tilde\Delta_\Omega$ is not symmetric.

\par {\bf 4.4. Quantum hypergroups associated with a twisting of the group
${\Bbb Z}^2_m\rtimes {\Bbb Z}_2$.}
Let  $G = {\Bbb Z}^2_m\rtimes {\Bbb Z}_2$, $m\geq3$, be a finite
group of order $2m^2$ with the following action $\alpha$ of ${\Bbb Z}_2 = \{\mbox{id}, s\}$
on $H = {\Bbb Z}_m^2 = \{(a, b)| a, b = 0, 1, \dots , m-1 \}$:  $\alpha_s (a, b) = (b, a)$.
The twisting of the group $G$ was constructed in \cite{Vai:twist} by using the 2-cocycle
$\omega$ on $\hat H = \hat{{\Bbb Z}}_m^2 \times \hat{{\Bbb Z}}_m^2$:
$$
\omega (\hat a, \hat b; \hat c, \hat d)
= \exp \left(\frac {2\pi i} m (\hat a\hat d - \hat b\hat c)\right).
$$
Let $\Omega$ be the counital 2-cocycle on ${\Bbb C}(G)\otimes {\Bbb C}(G)$ obtained by
the lifting construction:
$$
\Omega = \sum_{\hat H\times\hat H}
\exp\left(\frac {2\pi i} m (\hat a\hat d - \hat b\hat c)\right) P_{(\hat a, \hat b)}\otimes P_{(\hat c, \hat d)}.
$$
Then the twisted Kac algebra ${\cal A}_\Omega$ is not symmetric, for example,
$\Delta_\Omega (\lambda(s)) \neq \Sigma \Delta_\Omega (\lambda(s))$.

\par Define an automorphism $\gamma$ of the group $G$ as follows:
$\gamma (s) = s$ and $\gamma (a, b) = (ar, br)$, where $(a, b) \in H$
and $r^2 = 1$ modulo $m$. Then it is strightforward that
$(\gamma\otimes\gamma)(\Omega) = \Omega$. Thus $\gamma$ is an
automorphism of the twisted Kac algebra ${\cal A}_\Omega$ by
Proposition 3.6. Let $P$ be the conditional expectation on ${\cal A}_\Omega$
associated with the subgroup $\Gamma = \{\mbox{id}, \gamma\}$.
Then $B = P({\Bbb C}(G))$ is a quantum hypergroup.
Since $\gamma (s) = s$, the algebra $B$ is noncommutative.
The dimension of $B$ is equal to the number of $\Gamma$-orbits, hence
$\dim B = m^2 + p^2$, where $p$ is the number of solutions of the
equation $(r-1)k = 0$ modulo $m$.

\par If $r=m-1$, then the quantum hypergroup $B$ is non-symmetric.
Indeed, it is easy to see that $(\gamma\otimes\mbox{id})(\Omega) = \Omega^* =
(\mbox{id}\otimes\gamma)(\Omega)$ and $\Sigma\Omega = \Omega^*$. Then,
for any $a\in H$, we have
\begin{eqnarray*}
(P&\otimes & P)(\Delta_\Omega (\lambda(as)) - \Sigma\Delta_\Omega (\lambda(as)))\\
&=& \frac 1 4 (\mbox{id}\otimes\mbox{id}+\mbox{id}\otimes\gamma +\gamma\otimes\mbox{id}+\gamma\otimes\gamma)
(\Omega\lambda(as)\otimes\lambda(as)\Omega^* - \Omega^*\lambda(as)\otimes\lambda(as)\Omega)\\
&=& (\lambda(a) - \lambda(a^{-1}))\otimes (\lambda(a) - \lambda(a^{-1}))
\left(\Omega\lambda(s)\otimes\lambda(s)\Omega^* - \Omega^*\lambda(s)\otimes\lambda(s)\Omega\right),
\end{eqnarray*}
since $\Omega$ commutes with ${\Bbb C}(H)$. Thus $\hat \Delta_\Omega$ is
non-symmetric iff $(\lambda(a) - \lambda(a^{-1}))\otimes (\lambda(a) - \lambda(a^{-1}))
\left(\Omega^2\lambda(s)\otimes\lambda(s) - \lambda(s)\otimes\lambda(s)\Omega^2\right)\neq 0$.
This inequality can be obtained by strightforward calculations.

\bigskip

\noindent
{\bf Address:}

\bigskip

\noindent
        A. A. Kalyuzhnyi,
	Institute of Mathematics, Ukrainian
        National Academy of Sciences, ul. Tereshchinkivs'ka, 3, Kiev
        262601, Ukraine

\noindent
E-mail: kalyuz@imath.kiev.ua

\end{document}